\documentclass[draft,utf8,12pt,english,nothms]{aclart}

\usepackage{natbib}
\makeatletter \def\NAT@nmfmt#1{{\scshape\NAT@up#1}} \makeatother

\makeatletter
  \def\bibliography#1{%
  \if@filesw
    \immediate\write\@auxout{\string\bibdata{#1}}%
  \fi
  \expandafter\input{\bbl@main@language bst.tex}%
  \@input@{\jobname.bbl}}

\makeatother

\usepackage{smfhref}

\numberwithin{equation}{subsubsection}
\def\thesection{\arabic{section}}
\def\thesubsection{\thesection.\arabic{subsection}}
\let\subsection\Subsection

\input t2aenc.def
\def\BM{{\text{\upshape \fontencoding{T2A}\selectfont \char'301}}}

\usepackage[T1]{fontenc}
\usepackage{lmodern,fourier}
\usepackage{amssymb,graphicx}
\usepackage[matrix,arrow]{xy}

\theoremstyle{plain}

\theoremstyle{definition}

\let\mathsec\mathsf

\def\P{{\mathbf{P}}}
\def\C{{\mathbf C}}
\def\R{{\mathbf R}}
\def\Q{{\mathbf Q}}
\def\Z{{\mathbf Z}}

\def\Val{\operatorname{Val}}
\def\Gr{\operatorname{Gr}}
\def\GL{\operatorname{GL}}
\def\SL{\operatorname{SL}}
\def\eff{{\text{\upshape eff}}}
\def\ample{{\text{\upshape ample}}}
\def\AD{{\mathbb A}}
\def\Frob{\operatorname{Frob}}
\def\Tr{\operatorname{Tr}}

\title[Lectures on height zeta functions]{Lectures on height zeta functions: \\
At the confluence of algebraic geometry,
algebraic number theory, and analysis} 
\subtitle{French-Japanese Winter School, Miura (Japan), January 2008}
\author{Antoine Chambert-Loir}
\address{IRMAR, Universit\'e de Rennes~1 \\
Campus de Beaulieu \\ 35042 Rennes Cedex \\ France }
\email{antoine.chambert-loir@univ-rennes1.fr}

\def\C{{\mathbf C}}
\def\Q{{\mathbf Q}}
\def\R{{\mathbf R}}

\def\P{{\mathbf P}}
\def\Z{{\mathbf Z}}
\def\A{{\mathbf A}}

\def\abs#1{\left\lvert{#1}\right\rvert}
\def\Hom{\operatorname{Hom}}
\def\Re{\operatorname{Re}}

\def\norm#1{\left\lVert{#1}\right\rVert}

\def\hPic{\mathop{{\widehat {\operatorname{Pic}}}}}
\def\Pic{\operatorname{Pic}}

\def\gm{{\mathbf G_{\mathrm m}}}
\def\ga{{\mathbf G_{\mathrm a}}}

\let\bar\overline

\let\ra\rightarrow

\let\leq\leqslant \let\geq\geqslant
\let\phi\varphi\let\eps\varepsilon

\def\Res{\operatorname{Res}}
\def\Card{\operatorname{Card}}

\begin{document}
\selectlanguage{english}

\maketitle

\setcounter{tocdepth}2
\tableofcontents

\section{Introduction}
\subsection{Diophantine equations and geometry}

\subsubsection{Diophantine equations}
Broadly speaking, \emph{arithmetic} is the study of diophantine
equations, that is, systems of polynomial
equations with \emph{integral} coefficients,
with a special emphasis on  their solutions in rational integers.
Of course, there are numerous variants, the most obvious
ones allowing to consider coefficients and solutions
in the field of rational numbers, or in more general number
fields, or even in more general fields, \emph{e.g.},
finite fields.

The reader should be warned that,
in this generality, we are constrained 
by the \emph{undecidability theorem}
of~\cite{matiyasevich1971}:
there is no general method, that is no algorithm,
to decide whether or not any given polynomial system has
solutions in rational integers.
Any mathematician working on diophantine  equations
is therefore obliged to consider specific types
of diophantine equations, in the hope that such an undecidability
issues do not apply within the chosen families 
of equations.

\subsubsection{Enters geometry}
At first, 
one is tempted to sort equations according to the
degrees of the polynomials which intervene. However, this approach is
much too crude, and
during the \textsc{xx}th century, mathematicians
were led to realize that there are profound relations
between the given diophantine equation and the \emph{geometry} of 
its solutions in real or complex numbers.
This led to considerations of geometric invariants such as the \emph{genus}
of an algebraic curve, to the (essentially opposite) notions of a variety
of general type and a Fano variety, to the notion of a rational variety,
etc.

In this survey, we are interested in diophantine equations having
infinitely many solutions. A natural way to describe this infinite
set consists in sorting the solutions according to their size
(as integers) and in studying the asymptotic behaviour of the 
number of solutions of size smaller than a growing bound.

\subsubsection{The circle problem}
The classical \emph{circle problem} in analytic number theory
is to estimate the number of integer vectors $\mathbf x\in\Z^n$
such that $s(\mathbf x)\leq B$, when $B\ra\infty$ and $s(\cdot)$
is an appropriate notion of  a size of a vector in~$\R^n$.
When $s(\cdot)=\norm{\cdot}$ is a norm, for example the euclidean norm,
this amounts to counting the number of lattice points in a
ball with center~$0$ and of radius~$B$.

Since a ball is convex, the leading term is easily seen
to equal the volume of this ball; in other words,
\[ \Card \{\mathbf x\in\Z^n\,;\, \norm{\mathbf x}\leq B\}\sim B^n V_1, \]
where $V_1$ is the volume of the unit ball in~$\R^n$.
The study of the error term, however, is of a much more delicate nature.
For $n=2$, there is an easy $\mathrm O(B)$-bound
which only requires the Lipschitz property of the boundary
of the unit ball; when the norm is euclidean, one can prove
a $\mathrm O(B^{2/3})$ bound using the positivity of the curvature
of the boundary; however, the conjectured~$\mathrm O(B^{1/2})$-bound
remains open.

\subsection{Elementary preview of Manin's problem}
\subsubsection{Projective space}
As indicated above, the size 
of a solution of a diophantine equation can be thought of as a
measure of its complexity;
for a vector $\mathbf x\in\Z^n$, it 
is standard to define its size to be its euclidean norm. 

However, in algebra and geometry, we are often led to consider
rational functions and their poles, which inevitably bring us to ``infinity''.
The appropriate context to define and study the corresponding sizes, 
or as we shall now say, \emph{heights},
is that of \emph{projective geometry}.
Recall that for any field~$F$,
one defines the $n$-dimensional projective space $\P^n(F)$ on~$F$
to be the set of lines in the vector space~$F^{n+1}$.
In other words, this is the quotient set of 
the set $F^{n+1}\setminus\{0\}$ of nonzero vectors
modulo the action of homotheties:
a point in $\P^n(F)$ can be described 
by a nonzero family of $n+1$ \emph{homogeneous coordinates} $[x_0:\cdots:x_n]$,
while we are free to multiply these coordinates
by a common nonzero element of~$F$, still describing the same point.

To stick to the current terminology of algebraic geometry,
this set $\P^n(F)$ will be called the set of \emph{$F$-rational points
of the projective space~$\P^n$} and the latter will be
referred to as the \emph{scheme $\P^n$}.

\subsubsection{Heights}
Consider the case where $F=\Q$, the field of rational numbers.
Any point $\mathbf x\in\P^n(\Q)$ can be represented by
$n+1$ rational numbers, not all zero; however, if we multiply
these homogeneous coordinates by a common denominator,
we see that we may assume them to be integers; we then
may divide them by their greatest common divisor
and obtain a system of homogeneous coordinates $[x_0:\cdots:x_n]$
made of
$n+1$ coprime integers. At this point, only one choice
is left to us, namely multiplying this system by~$-1$.

Consequently, we may define the \emph{exponential height}
of~$\mathbf x$ as
$ H(\mathbf x) = \max(\abs{x_0},\dots,\abs{x_n}) $
and its \emph{logarithmic height} as $h(\mathbf x)=\log H(\mathbf x)$.
(Observe the notation, popularized  
by Serge \textsc{Lang}: small ``h'' for logarithmic height,
capital ``H'' for exponential height.)

\subsubsection{The theorems of {Northcott} and {Schanuel}}
\cite{northcott1950} made a fundamental, albeit trivial,
observation: for any real number~$B$, there
are only finitely many points $\mathbf x\in\P^n(\Q)$
such that $H(\mathbf x)\leq B$.
Indeed, this amounts to saying that there are only
finitely many systems of coprime integers~$[x_0:\cdots:x_n]$
such that $\abs {x_i}\leq B$ for all~$i$,
considered up to multiplication by~$\pm1$.
Concretely, there are at most $(2B+1)^{n+1}$ such systems
of integers, let it be coprime, and modulo~$\pm1$!

More precisely, \cite{schanuel79} proved that the number~$N(B)$
of such points satisfies
\[ N(B)  \sim \frac{2^n}{\zeta(n+1)} B^{n+1} \]
when $B\ra\infty$, where $\zeta(n+1)$ is Riemann's zeta function
evaluated at~$n+1$.
When, say, $n=1$, this can be interpreted as follows: for
$B\ra\infty$, the probability that two integers from~$[1,B]$ are coprime
tends to $1/\zeta(2)=6/\pi^2$, a result originally due to \textsc{Dirichlet};
see Theorem~332 in the classic book of~\cite{hardy-wright1979}.

The remarkable appearance of this notorious function reveals
that something profound is happening before our eyes,
something that would certainly appeal to any arithmetician.

\subsubsection{The problem of Batyrev--Manin}
More generally, we want to consider subsets of the projective
space defined by polynomial equations.
Evaluating a polynomial
in~$n+1$ variables at a system of homogeneous coordinates
of a point $\mathbf x\in\P^n(F)$ does not lead
to a well-defined number, since the result changes very much 
if we multiply the coordinates by a common nonzero factor~$\lambda$.
However, if our polynomial is assumed to be homogeneous
of some degree, the result is multiplied by~$\lambda^d$;
at least we can consistently say whether the result is zero or not.

Let thus $X\subset\P^n$ be a \emph{closed subscheme} of~$\P^n$,
namely the set of common zeroes of a family of homogeneous
polynomials in~$n+1$ variables. In other words, for any field~$F$,
$X(F)$ is the subset of~$\P^n(F)$ consisting of points
satisfying this family of equations.
In fact, we will also consider the case
of locally closed subschemes of~$\P^n$, that is, subschemes 
of the form $X\setminus Y$, where $Y$
is a closed subscheme of~$\P^n$ contained in a closed subscheme~$X$.

We now define $N_X(B)$ to be the number of points~$\mathbf x\in X(\Q)$
such that $H(\mathbf x)\leq B$. 
The problem of~\cite{batyrev-m90} is to understand the
asymptotic behaviour of~$N_X(B)$, for $B\ra\infty$. In particular,
one is interested in the two numbers:
\[ \beta_X^- = \liminf_{B\ra\infty}  \frac{\log N_X(B)}{\log B},\qquad
   \beta_X^+ = \limsup _{B\ra\infty}\frac{\log N_X(B)}{\log B}.\]
Observe that $\beta_X^-\leq\beta_X^+$
and that, when $X=\P^n$, $\beta^-_X=\beta^+_X=n+1$.

\subsubsection{Heath-Brown's conjecture}
However, some problems in analytic number theory require uniform
upper bounds, e.g., bounds depending only on the degree and dimension of 
$X$.
\cite{heath-brown2002} conjectured
that for any real number~$\eps>0$ and any integer~$d>0$,
there exists a constant $C(d,\eps)$ such that $N_X(B)\leq C(d,\eps) B^{\dim X+\eps}$,
for any integral closed subscheme $X$ of~$\P^n$ of degree~$d$ and 
dimension~$\dim X$. (Degree and dimensions are defined
in such a way that when $H$ is a general subspace of codimension~$\dim X$,
$(X\cap H)(\C)$ is a finite set of cardinality~$\deg X$, counted with
multiplicities.  Here,
``integral'' means that $X$ is not a nontrivial union of two subschemes.)

As an example for the conjecture, let us consider the polynomial
$f=x_0x_2-x_1x_3$ and the subscheme $X=V(f)$ in~$\P^3$ defined by~$f$.
Considering lines on~$X$, it may be seen that $N_X(B)\gg B^2$.


In fact, \cite{browning-hb-s2006} established
that this conjecture is equivalent to the same statement 
with $X$ assumed to be a hypersurface
defined by a single homogeneous polynomial of degree~$d$
and proved the case $d\geq 6$. The strategy  of the proof 
comes from~\cite{heath-brown2002} (from which
the case $d=2$ also follows) and uses
previous work of~\cite{bombieri-pila1989} about
integral points on plane curves.
Consequently, the only open cases are those for which $3\leq d\leq 5$
and are the object of active works: to quote only two of them,
the case of smooth hypersurfaces has been completed
in~\cite{browning-hb2006},
while \cite{salberger2007} treats the case where
$d\geq 4$ and $X$ contains only finitely many linear
spaces of codimension~$1$.
As this is not the main topic of this survey,
we refer to these articles for more details and references.

\subsubsection{Geometric parameters}
Let us discuss the parameters at our disposal
to describe the counting function~$N_X(B)$
and the exponents~$\beta_X^-$ and~$\beta_X^+$.

The most important  of them come from geometry. First of all, there is the dimension,
as was already seen in \textsc{Heath-Brown}'s conjecture.

The dimension appears also in \textsc{Schanuel}'s theorem (precisely,
the dimension plus one), but this is not a good interpretation
of the exponent. Namely, one of the insights of~\cite{batyrev-m90}
was the interpretation of~$n+1$ as describing the location
of the anticanonical divisor of~$\P^n$ with respect to the cone
of  effective divisors. In effect, the $n$-differential form
\[ d(x_1/x_0) \wedge d(x_2/x_0)\wedge\dots\wedge d(x_n/x_0) \]
has a pole of multiplicity~$n+1$ along the hyperplane 
defined by the equation~$x_0=0$.
Note that the rank of the N\'eron--Severi group 
will also  enter the final picture --- 
this is the group of classes of divisors where we identiy
two divisors if they give the same degree to any curve drawn on
the variety.

Properties of a more arithmetic nature intervene as well:
the classical Hasse principle and weak approximation,
but also their refinements within the framework of the theory of 
\emph{Brauer-Manin obstructions}. This explains the appearance of orders of
various Galois cohomology groups in asymptotic formulas.

Were we to believe that the result only depends on
the ``abstract'' scheme~$X$, we would rapidly find ourselves
in contradictions. Indeed, there are many $\P^1$ 
that can be viewed in~$\P^2$. On lines, the 
number of points  of bounded height has the same order of magnitude,
but the leading coefficient can change;
for example, for lines with equations, say $x_0=0$, or $x_0+x_1+x_2=0$,
one finds $\beta^{\pm}=2$,
but the leading coefficient in~$N_X(B)$ is $2/\zeta(2)$ in the first case, and  
only one half of it in the second case.
There are also ``non-linear lines'' in~$\P^2$, namely 
curves which can be parametrized by rational functions of one variable.
The equation $x_2x_0-x_1^2$ defines such a curve~$X'$
for which $\beta_{X'}^{\pm}=1$, reflecting
its non-linearity.
In other words, even if we want to think in terms of an abstract variety~$X$,
the obtained formulas 
force us to take into account the \emph{embedding} of~$X$ in a
projective space~$\P^n$.

There are also more subtle geometric and arithmetic caveats,
related not only to~$X$, but to its closed subschemes;
we will consider these later.

\subsubsection{The height zeta function}
Seeking tools and applications of analysis,
we introduce the \emph{height zeta function}.
This is nothing other than the generating Dirichlet series
$Z_X(s) = \sum_{\mathbf x\in X(\Q)} H(\mathbf x)^{-s}$,
where $s$ is a complex parameter.

This series converges for $\Re(s)$ large enough; for example,
\textsc{Schanuel}'s theorem implies that it converges for $\Re(s)>n+1$
and that it defines a holomorphic function on
the corresponding half-plane.

A first, elementary, analytic invariant related to~$Z_X$ is its \emph{abscissa
of convergence}~$\beta_X$.
Since this Dirichlet series has positive coefficients,
a result of \textsc{Landau} implies the inequality
\[ \beta_X^-\leq \beta_X \leq \beta_X^+. \]
This is a very crude form of Abelian/Tauberian theorem,
which we will use later on, 
once we have obtained precise analytic information
about~$Z_X$.
Our goal will be to establish an analytic continuation
for~$Z_X$ to a meromorphic function on some larger domain of
the complex plane. Detailed information on the poles of
this continuation will then imply a precise asymptotic
expansion for~$N_X(B)$.

\section{Heights}
\subsection{Heights over number fields}
\subsubsection{Absolute values on the field of rational numbers}

Let us recall that an absolute value~$\abs\cdot$ on a field~$F$
is a function from~$F$ to~$\R_+$ satisfying the following properties:
\begin{itemize}
\item $\abs a=0$ if and only if $a=0$ (non-degeneracy);
\item $\abs {a+b}\leq\abs{a}+\abs b$ for any $a,b\in F$ (triangular inequality);
\item $\abs{ab}=\abs a\abs b$ for any $a,b\in F$ (multiplicativity).
\end{itemize}

Of course, the usual modulus on the field of complex numbers
satisfies these properties, hence is an absolute value on~$\C$.
It induces an absolute value on any of its subfields, in particular
on~$\Q$. It is called archimedean since for any $a\in\Q$
such that $\abs a\neq 0$, and any~$T>0$,
there exists an integer~$n$ such that $\abs{na}>T$. We write~$\abs\cdot_\infty$
for this absolute value.

In fact, the field of rational numbers  possesses many other absolute
values, namely the $p$-adic absolute value, where $p$ is any prime number.
It is defined as follows: Any nonzero rational number~$a$
can be written as $p^m u/v$, where $u$ and~$v$ are integers
not divisible by~$p$ and~$m$ is a rational integer; 
the integer $m$ depends only on~$a$ and we define $\abs a_p=p^{-m}$;
we also set $\abs 0_p=0$.
Using uniqueness of factorization of integers into prime numbers,
it is an exercise to prove that $\abs\cdot_p$ is an absolute
value on~$\Q$. In fact, not only does the triangular inequality 
hold, but a stronger form is actually true: the \emph{ultrametric}
inequality:
\begin{itemize}
\item $\abs{a+b}_p\leq\max(\abs a_p,\abs b_p)$, for any $a,b\in\Q$.
\end{itemize}

\subsubsection{Ostrowski's theorem}
The preceding list gives us essentially all absolute values on~$\Q$.
Indeed, let $\abs{\cdot}$ be an absolute value on~$\Q$.
By \emph{\textsc{Ostrowski}'s theorem}, $\abs\cdot$ is one of the following (mutually
exclusive) absolute values:
\begin{itemize}
\item the trivial absolute value $\abs\cdot_0$,
defined by $\abs a_0=1$ if $a\neq0$ and $\abs 0_0=0$;
\item the standard archimedean absolute value $\abs\cdot_\infty$
and its powers $\abs\cdot_\infty^s$ for $0<s\leq 1$;
\item the $p$-adic absolute value $\abs\cdot_p$
and its powers $\abs\cdot_p^s$ for $0<s<\infty$
and some prime number~$p$.
\end{itemize}

\subsubsection{Topologies and completions}
To each of these families corresponds some distance on~$\Q$
(trivial, archimedean, $p$-adic), given by $d(a,b)=\abs{a-b}$,
hence some topology. The trivial absolute value defines
the discrete topology on~$\Q$, the archimedean absolute values
the usual topology of~$\Q$, as a subspace of the real numbers,
and the $p$-adic absolute values the so-called $p$-adic topology.

The standard process of Cauchy sequences now constructs
for any of these absolute values
a complete field in which~$\Q$ is dense and to which the
absolute value uniquely extends.
The obtained field is~$\Q_\infty=\R$ when the absolute value
is the archimedean one, and is written~$\Q_p$ when the
absolute value is $p$-adic.

\subsubsection{Product formula}
A remarkable equality relates all of these absolute values, namely:
for any nonzero rational number $a$, $\prod_{p\leq\infty} \abs a_p=1$.
This is called the \emph{product formula}.
(We have written $p\leq\infty$ to indicate that the set of indices~$p$
ranges over the set of prime numbers to which we adjoin the symbol~$\infty$.)
Let us decompose $a=\pm \prod_p p^{n_p}$ as a product of
a sign and of (positive or negative) prime powers. For any prime number~$p$,
we have $\abs a_p=p^{-n_p}$, while $\abs a_\infty=\prod_p p^{n_p}$;
the product formula follows at once from these formulae.

From now on, we shall ignore the trivial absolute value.

\subsubsection{Number fields}
Let $F$ be a number field, that is, a finite extension of~$\Q$.
Let $\abs\cdot_v$ be a (non-trivial)
absolute value on~$F$. Its restriction to~$\Q$
is an absolute value on~$\Q$, hence is given  by a power of
the $p$-adic absolute value (for $p\leq\infty$), as in
\textsc{Ostrowski}'s list. The completion furnishes a complete field~$F_v$
which is again a finite extension of~$\Q_p$; in particular, we have
 a norm map $\mathrm N\colon F_v\ra\Q_p$. We say that $v$ 
is \emph{normalized} if $\abs a_v=\abs {N(a)}_p$.

We write $\Val(F)$ for the set of (non-trivial) normalized
valuations on~$F$. As in the case of~$\Q$, the product
formula holds: for any nonzero $a\in F$,
\[ \prod_{v\in\Val(F)} \abs a_v = 1. \]

\subsubsection{Heights on the projective space}
Using this machinery from algebraic number theory,
we may extend the definition of the height of a point in  $\P^n(\Q)$
to points in~$\P^n(F)$.

Let $\mathbf x$ be a point of~$\P^n(F)$, given by a system
of homogeneous coordinates~$[x_0:\dots:x_n]$ in~$F$, not all zero.
We may define the (exponential) height of~$\mathbf x$ as
\[ H_F(\mathbf x) = \prod_{v\in\Val(F)} \max(\abs{x_0}_v,\abs{x_1}_v,
\dots, \abs{x_n}_v ) \]
since the right hand side does not depend on the choice of a specific
system of homogeneous coordinates. Indeed,
if we replace~$x_i$ by~$ax_i$, for some nonzero element of~$F$,
the right hand side gets multiplied by
\[  \prod_{v\in\Val(F)}  \abs{a}_v=1 .\]

For $F=\Q$, this definition coincides with the one previously 
given. Indeed, we may assume that the coordinates~$x_i$
are coprime integers. Then, for any prime number~$p$,
$\abs{x_i}_p\leq 1$ (since the $x_i$s are integers), and one
of them is actually equal to~$1$ (since they are not all divisible by~$p$).
Consequently, the $p$-adic factor is equal to~$1$ and
\[ H_\Q(\mathbf x) = \max(\abs{x_0}_\infty,\dots,\abs{x_n}_\infty) = H(\mathbf x). \]

Using a bit more of algebraic number theory, one can show that
for any finite extension~$F'$ of~$F$, and any $\mathbf x\in\P^n(F)$,
\[ H_{F'}(\mathbf x)= H_F(\mathbf x)^d, \qquad\text{where $d=[F':F]$.}\]

\subsection{Line bundles on varieties}

\subsubsection{The notion of a line bundle}
One of the main differences between classical
and modern algebraic geometry is the change of emphasis
from ``subvarieties of a projective space'' to 
``varieties which can be embed in a projective space''.
Given an abstract variety~$X$, the data of an embedding is 
essentially described by line bundles, a notion we now have
to explain.

Let $X$ be an algebraic variety over some field~$F$. 
A \emph{line bundle} $L$ on~$X$
can be thought of as a ``family'' of lines~$L(x)$ 
(its \emph{fibres}) parametrized by the points~$x$ of $X$. 
There is actually a subtle point
that not only rational points $x\in X(F)$ must be considered,
but also points  in $X(F')$, for all extensions~$F'$ (finite or not)
of the field~$F$. (We will try to hide these kind of complications.)

A line bundle $L$ has \emph{sections}: over an open subset~$U$
of~$X$, such a section~$s$ induces maps $x\mapsto s(x)$ for $x\in U$,
where $s(x)\in L(x)$ for any~$x$. A section can be multiplied
by a regular function: if $s$ is a section of~$L$ on~$U$
and $f$ is a regular function on~$U$, then there is a section $fs$
corresponding to the assignement $x\mapsto f(x)s(x)$.
At the end, the set $\Gamma(U,L)$ of sections of~$L$ 
on~$U$, also written $L(U)$, has a natural structure
of a \emph{module} on the ring~$\mathscr O_X(U)$ of regular
functions on~$U$.

Sections can also be glued together: if $s$ and~$s'$ are sections of~$L$
over open subsets~$U$ and~$U'$ which coincide on~$U\cap U'$,
then there is a unique section~$t$ on~$U\cap U'$ which
induces~$s$ on~$U$ and~$s'$ on~$U'$.

What ties all of these lines together is that for every point~$x\in X$,
there is an open neighbourhood~$U$ of~$x$, a \emph{frame} $\eps_U$
which is a section of~$L$ on~$U$ such that $\eps_U(x)\neq 0$
for all $x\in U$ and such that any other section $s$ of~$L$ on~$U$
can be uniquely written as $f \eps_U$, where $f$
is a regular function on~$U$.

\subsubsection{The canonical line bundle}
The sections of the \emph{trivial line bundle} $\mathscr O_X$
on an open set~$U$ are the ring
$\mathscr O_X(U)$ of regular functions on~$U$.

If $X$ is smooth and everywhere $n$-dimensional, 
it possesses a canonical  line bundle~$\omega_X$,
defined in such a way that its sections on an open subset~$U$
are precisely the module of $n$-differential forms,
defined algebraically as the $n$-th exterior product
of the $\mathscr O_X(U)$-module of Kähler differentials
$\Omega^1_{\mathscr O_X(U)/F}$.

\subsubsection{The Picard group} 
There is a natural notion of morphism of line bundles.
If $f\colon L\ra M$ is such a morphism, it assigns to any section
$s\in \Gamma(U,L)$ a section $f(s)\in\Gamma(U,M)$
so that the map $s\mapsto f(s)$ is a morphism of $\mathscr O_X(U)$-modules
(\emph{i.e.}, additive and compatible with multiplication by regular
functions).

An isomorphism is a morphism~$f$ for which there is an ``inverse morphism''
$g\colon M\ra L$ such that $g\circ f$ and $f\circ g$ are the identical
morphisms of~$L$ and~$M$ respectively.

From two line bundles~$L$ and~$M$ on~$X$, one can construct  a third
one, denoted $L\otimes M$, in such a way that
if $\eps$ and~$\eta$ are frames of~$L$ and~$M$ on an open set~$U$,
then $\eps\otimes\eta$ is a frame of~$L\otimes M$ on~$U$,
with the obvious compatibilities suggested by the tensor
product notation, namely
\[ (f\eps)\otimes\eta= \eps\otimes (f\eta)=f (\eps\otimes\eta), \]
for any regular function on~$U$.

Any line bundle~$L$ has an ``inverse'' $L^{-1}$ for the tensor product;
if $\eps$ is a frame of~$L$ on~$U$, then a frame of~$L^{-1}$ on~$U$
is $\eps^{-1}$, again with the obvious compatibilities
\[ (f\eps)^{-1}= f^{-1} \eps^{-1} \]
for any nonvanishing regular function~$f$ on~$U$.

These two laws are in fact compatible with the notion of isomorphism.
The set of isomorphism classes of line bundles on~$X$ 
is an Abelian group, which is called the \emph{Picard group}
and denoted $\Pic(X)$.

\subsubsection{Functoriality}
Finally, if $u\colon X\ra Y$ is a morphism of algebraic varieties
and $L$ is a line bundle on~$Y$, there is a line bundle $u^*L$
on~$X$ defined in such a way that the fibre of~$u^*L$ at
a point~$x\in X$ is the line $L(u(x))$; similarly,
if $\eps$ is a frame of~$L$ on an open set~$U$ of~$Y$,
then $u^*\eps$ is a frame of~$u^*L$ on~$u^{-1}(U)$.

At the level of isomorphism classes, this induces
a map $u^*\colon \Pic(Y)\ra\Pic(X)$ which is a morphism
of Abelian groups.

\subsection{Line bundles and embeddings}
\subsubsection{Line bundles on~$\P^n$}
As a scheme over a field~$F$, $\P^n$ parametrizes hyperplanes
of the fixed vector space~$V=F^{n+1}$. (In the introduction,
we identified a point $\mathbf x\in\P^n(F)$ with the \emph{line}
in~$V$ generated by any system of homogeneous coordinates;
this switch of point of view can be restored by duality.) 
In particular, to a point $x\in \P^n(F)$ corresponds a hyperplane~$H(x)$
in~$V$. Considering the quotient vectorspace $V/H(x)$,
we hence get a line $L(x)$.

These lines $L(x)$ form a line bundle on~$\P^n$.
This line bundle admits sections $s_0,\ldots,s_n$ defined on~$\P^n$
which correspond to the homogeneous coordinates.
Observe however that for $x\in \P^n(F)$, $s_i(x)$ is not  a number,
but a member of some line. However, given any generator~$\eps(x)$,
of this line, there are elements~$x_i\in F$ such that 
$s_i(x)=x_i\eps(x)$ and the family $[x_0:\cdots:x_n]$ 
gives homogeneous coordinates for~$x$. 

It is useful to remember
the relation $x_j s_i(x)=x_i s_j(x)$, valid for any $x\in\P^n(F)$
with homogeneous coordinates~$[x_0:\dots:x_n]$ and any couple~$(i,j)$
of indices. It will also be necessary to observe
that for any~$x\in\P^n$, at least one of the~$s_i(x)$ is nonzero
--- this is a reformulation of the fact that the homogeneous coordinates
of a point are not all zero.

In traditional notation, this line bundle is denoted~$\mathscr O_{\P^n}(1)$.
Its class in the Picard group of~$\P^n$ is a generator of this group,
which is isomorphic to~$\Z$. Similarly, the traditional notation
for the line bundle corresponding to an integer~$a\in\Z$ is 
$\mathscr O_{\P^n}(a)$.

The canonical line bundle of~$\P^n$ is isomorphic to~$\mathscr O_{\P^n}(-n-1)$.
In fact, in the open subset~$U_i$ of~$\P^n$ where the homogeneous coordinate $x_i$ is nonzero, $x_j/x_i$ defines a regular function and 
one has a differential form
\[ \omega_i = d(x_0/x_i) \wedge d(x_1/x_i) \wedge\dots
 \wedge \widehat{d(x_i/x_i)} \wedge\dots\wedge d(x_n/x_i), \]
where the hat indicates that one omits the corresponding factor.

On the intersections $U_i\cap U_j$, one can check
that the expressions
$s_i(x)^{\otimes n+1}\omega_i$ and $s_j(x)^{\otimes n+1}\omega_j$
identify the one to the other if one uses the relation $x_i s_j(x)=x_j s_i(x)$
that we observed. This implies that these expressions
can be glued together as a section of $\mathscr O_{\P^n}(n+1)\otimes \omega_{\P^n}$ which vanishes nowhere, thereby establishing the announced isomorphism.

\subsubsection{Morphisms to a projective space}
Now we translate into the language of line bundles
the geometric data induced  by a morphism $f$ from a variety~$X$
to a projective space~$\P^n$.
As we have seen, $\P^n$ is given with its line bundle~$\mathscr O_{\P^n}(1)$
and its sections $s_0,\dots,s_n$ which do not vanish simultaneously. 
By functoriality, the morphism~$f$ furnishes a line bundle
$f^*\mathscr O_{\P^n}(1)$ on~$X$ together with~$(n+1)$ sections
$f^*s_0,\dots,f^*s_n$ which, again, do not vanish simultaneously.

A fundamental fact in projective geometry is that this assignment 
defines a bijection between:
\begin{itemize}
\item the set of morphisms $f\colon X\ra\P^n$;
\item the set of data $(L,u_0,\dots,u_n)$ consisting of a line 
bundle~$L$ on~$X$
together with $n+1$ sections $u_0,\dots,u_n$ which do not vanish simultaneously.
\end{itemize}
We have only described one direction of this bijection,
the other can be explained as follows.
For $i\in\{0,\dots,n\}$, let $U_i$ be the open subset of~$X$
where $u_i\neq 0$; by assumption, $U_0\cup\dots\cup U_n=X$.
On~$U_i$, $u_i$ is a frame of~$L$ and there are regular functions
$f_{i0},\dots,f_{in}$ on~$U_i$ such that, on that open set,
$u_j=f_{ij} u_i$. This allows to define a morphism
of algebraic varieties $f_i\colon U_i\ra\P^n$ by the formula
$x\mapsto [f_{i0}(x):\dots:f_{in}(x)]$. It is easy to check
that for any couple~$(i,j)$,
the morphisms $f_i$ and~$f_j$ agree on~$U_i\cap U_j$,
hence define a morphism $f\colon X\ra\P^n$.

A simple example is given by the line bundle~$\mathscr O_{\P^1}(d)$
on~$\P^1$
and the sections $s_0^{\otimes i_0}\otimes s_1^{\otimes i_1}$, for $i_0+i_1=d$.
It corresponds to the Veronese embedding of~$\P^1$ in~$\P^d$
defined by $[x:y]\mapsto [x^d:x^{d-1}y:\dots:xy^{d-1}:y^d]$.

\subsubsection{Cones in the Picard group}
Let $X$ be an algebraic variety and $\Pic(X)_\R$ 
the real vector space obtained
by tensoring the group~$\Pic(X)$ with the field of real numbers.
Although we shall only study
cases where it is finite-dimensional, this vector space
may very well be infinite-dimensional.

The vector space  $\Pic(X)_\R$ contains two distinguished cones.
The first, $\Lambda_\eff$, is called the \emph{effective cone}, 
it is generated by line bundles which have nonzero sections over~$X$.
The second, $\Lambda_\ample$, is called the \emph{ample cone} 
and it is generated by line bundles of the form~$f^*\mathscr O_{\P^n}(1)$,
where $f$ is an \emph{embedding} of~$X$ into a projective space~$\P^n$.
(Such line bundles are called \emph{very ample}, a line bundle
which has a very ample power is called \emph{ample}.)

One has $\Lambda_\ample\subset\Lambda_\eff$, because a very ample line bundle
has nonzero sections, but the inclusion is generally strict.

For $X=\P^n$ itself, one has $\Pic(X)\simeq \Z$, hence $\Pic(X)_\R\simeq \R$,
the line bundle $\mathscr O_{\P^n}(1)$ corresponding to the number~$1$,
and both cones are equal to~$\R_+$ under this identification.

\subsection{Metrized line bundles}

\subsubsection{Rewriting the formula for the height}
Since an abstract variety may be embedded in many ways in a projective space,
there are as many possible definitions for a  height function on it.
To be able to explain what happens, we first will rewrite
the definition of the height of a point $\mathbf x\in\P^n(F)$
(where $F$ is a number field)
using line bundles.

Recall that $\Val(F)$ is the set of normalized absolute values on
the number field~$F$.
Let $\mathbf x\in\P^n(F)$ and let $[x_0:\dots:x_n]$ be a system
of homogeneous coordinates for~$\mathbf x$.
If $x_i\neq 0$, we may write, for any absolute value $v\in\Val(F)$,
\[ \max(\abs{x_0}_v,\dots,\abs{x_n}_v)
 = \left( \frac{ \abs{x_i}_v }{\max(\abs{x_0}_v,\ldots,\abs{x_n}_v)} \right)^{-1}
\abs{x_i}_v.\]
Recalling that~$s_i$ is the section of~$\mathscr O_{\P^n}(1)$
corresponding to the $i$th homogeneous coordinate, we define
\[ 
 \norm{s_i(\mathbf x)}_v =  \frac{ \abs{x_i}_v }{\max(\abs{x_0}_v,\ldots,\abs{x_n}_v)}.
\]
We observe that the right hand side does not depend on the choice
of homogeneous coordinates. 
Moreover, the product formula implies that
\[ 
H_F(\mathbf x) = \prod_{v\in\Val(F)} \norm{s_i(\mathbf x)}_v^{-1}
\abs{x_i}_v =  \prod_{v\in\Val(F)} \norm{s_i(\mathbf x)}_v^{-1}.\]
In other words, we have given a formula for $H_F(\mathbf x)$
as a product over~$\Val(F)$ where each factor is well defined,
independently of any choice of homogeneous coordinates.

The reader should not rush to conclusions: so far
we have only exchanged the indeterminacy of homogeneous
coordinates with the choice of a specific index~$i$,
more precisely of a specific section~$s_i$.

\subsubsection{Metrized line bundles: an example}
The notation introduced for $\norm{s_i(\mathbf x)}$ suggests
that it is the \emph{$v$-adic norm} of the vector $s_i(\mathbf x)$
in the fibre of~$\mathscr O(1)$ at~$\mathbf x$.
A $v$-adic norm on a $F$-vector space~$E$ is a map
$\norm{\cdot}_v\colon E\ra\R_+$ satisfying the following relations,
analogous to those that the $v$-adic absolute value possesses:
\begin{itemize}
\item $\norm{e}_v=0$ if and only if $e=0$ (non-degeneracy);
\item $\norm{e+e'}_v\leq\norm{e}_v+\norm{e'}_v$ for any $e,e'\in E$ (triangular inequality);
\item $\norm{ae}_v=\abs{a}_v \norm{e}_v$ for any $e\in E$ and any $a\in F$
(homogeneity).
\end{itemize}
 
In our case, the vector space is the line $E=\mathscr O(1)(\mathbf x)$
and the norm of a single nonzero vector in~$E$ determines
the norm of any other.
Consequently, the formula given for $\norm{s_i(\mathbf x)}_v$ (which
is a positive real number) uniquely extends to a norm on~$E$.
The formula $x_j s_i(\mathbf x)=x_is_j(\mathbf x)$
implies that the given norm does not depend on the initial choice
of an index~$i$ such that $s_i(\mathbf x)\neq 0$.

Moreover, when the point~$\mathbf x$ varies, the so-defined
norms vary continuously in the sense that the norm of a section~$s$
on a Zariski open set~$U$
extends to  a continuous function from~$U(F_V)$ 
(endowed with the $v$-adic topology) to~$\R_+$. 
It suffices to check this fact on open sets which
cover~$\P^n$ and over which $\mathscr O(1)$ has frames;
on the set where $x_i\neq 0$, $s_i$ is such a frame and the claim
follows by observing that the given
formula is continuous on~$U_i(F_v)$.

\subsubsection{Metrized line bundles: definition}
We now extend the previous construction to a general definition.
Let $X$ be an algebraic variety over a number field~$F$
and let $L$ be a line bundle on~$X$.
A \emph{$v$-adic metric} on~$L$ is the data of $v$-adic norms
on the $F_v$-lines $L(x)$, when $x\in X(F_v)$, which
vary continuously with the point~$x$. This assertion
means that for any open set~$U$ and any section~$s$ of~$L$ on~$U$,
the function $U(F_v)\ra\R_+$ given by $x\mapsto \norm{s(x)}_v$
is continuous. Since the absolute value of a regular function
is continuous, it suffices to check this fact for frames
whose open sets of definition cover~$X$.
A particular type of metrics is important; we call them \emph{smooth metrics}.
These are the metrics such that the norm of a local frame
is $\mathscr C^\infty$ if $v$ is archimedean, and locally constant
if $v$ is ultrametric.

The construction we have given in the preceding Section
therefore defines a $v$-adic metric on the line bundle~$\mathscr O_{\P^n}(1)$
on~$\P^n$. This metric is smooth if $v$ is ultrametric,
but not if $v$ is archimedean, because the function $(x,y)\mapsto \max(\abs x,\abs y)$ from~$\R^2$ to~$\R_+$ is not smooth.
A variant of the construction furnishes a smooth metric in that case,
namely the Fubini-Study metric, defined by
replacing $\max(\abs{x_0}_v,\ldots,\abs{x_n}_v)$ by 
$(\abs{x_0}_v^2+\dots+\abs{x_n}_v^2)^{1/2}$.

\subsubsection{Adelic metrics}
The formula for the height features all normalized absolute values of
the number field~$F$. We thus define an \emph{adelic metric} on a line bundle~$L$
to be a family of $v$-adic metrics on~$L$, for all $v\in\Val(F)$.
However, to be able to define a height using such data,
we need to impose a compatibility condition  of ``adelic type''
on all these metrics.

Let us describe this condition.
Let $U$ be an affine open subset of~$X$ and~$\eps$ be a frame of~$L$ on~$U$.
Let us represent~$U$ as a subvariety of some affine space~$\A^N$
defined by polynomials with coefficients in~$F$.

For any ultrametric absolute value~$v\in\Val(F)$, the subset
$\mathfrak o_v\subset F_v$ consisting of elements $a\in F_v$
such that $\abs{a}_v\leq 1$ is therefore a subring of~$F_v$.
We may thus consider the subset $U(\mathfrak o_v)=U(F_v)\cap \mathfrak o_v^N$ 
of $U(F_v)$.
Although it depends
on the specific choice of a representation of~$U$ as a subvariety
of~$\A^N$, one can prove that two representations will
define the same subsets $U(\mathfrak o_v)$ up to finitely many exceptions
in~$\Val(F)$. As a consequence, for any $x\in U(F)$,
one has $x\in U(\mathfrak o_v)$ up to finitely many exceptions~$v$
(choose a representation where~$x$ has coordinates~$(0,\dots,0)$).

The \emph{adelic compatibility condition} can now be expressed
by requiring  that for all $v\in\Val(F)$,
up to finitely many exceptions, $\norm{\eps(x)}_v=1$ for any $x\in U(\mathfrak o_v)$.

\subsubsection{Heights for adelically metrized line bundles}
Let $X$ be a variety over a number field $F$ and $L$ a line bundle
on~$X$ with an adelic metric.
For any $x\in X(F)$ and any frame $s$ on a neighbourhood~$U$ of~$x$,
let us define
\[ H_{L,s}(x) = \prod_{v\in\Val(F)} \norm{s(x)}_v^{-1}. \]
By definition of an adelic metric, almost all of the terms
are equal to~$1$.
If $t$ is any nonzero element of~$L(x)$, there exists
$a\in F^*$ such that $t=as(x)$; then $\norm{t}_v=\abs{a}_v\norm{s(x)}_v$
for any $v\in\Val(F)$. In particular, $\norm t_v=1$ for almost all~$v\in\Val(F)$
and 
\begin{align*} \prod_{v\in\Val(F)} \norm{t}_v^{-1}
& = \prod_{v\in\Val(F)} \norm{as(x)}_v^{-1} 
 =\prod_{v\in\Val(F)} \abs{a}_v^{-1} \norm{s(x)}_v^{-1} \\
& =\left(\prod_{v\in\Val(F)} \abs{a}_v^{-1}\right)\left( \prod_{v\in\Val(F)}\norm{s(x)}_v^{-1}\right) 
 = \prod_{v\in\Val(F)}\norm{s(x)}_v^{-1}= H_{L,s}(x)\end{align*}
where we used the product formula to establish the penultimate equality.

As a consequence, one can use any nonzero element
of~$L(x)$ in the formula $H_{L,s}(x)$; the resulting
product does not depend on the choice of~$s$. We write it $H_L(x)$.

For $L=\mathscr O_{\P^n}(1)$ with the adelic metric constructed
previously, we recover the first definition of the height.
As we shall see, other choices of adelic metrics furnish
the same function, up to a multiplicative factor
which is bounded above and below.
However, they can be of significative arithmetical interest.

\subsubsection{Properties of metrized line bundles}
There is a natural notion of tensor product of adelically
metrized line bundles, for which the $v$-adic norm
of a tensor product $e\otimes e'$ (for $e\in L(x)$ and $e'\in L'(x)$)
is equal to $\norm e_v \norm{e'}_v$. 
With that definition,
one has
\[ H_{L\otimes L'}(x)=H_L(x) H_{L'}(x). \]

An isometry of adelically metrized line bundles is an isomorphism
which preserves the metrics.  The set of isometry classes
of metrized line bundles form a group, called the \emph{Arakelov-Picard group},
and denoted $\hPic(X)$.
The above formula implies that the function from~$X(F)\times \hPic(X)$
to $\R_+^*$ is linear  in the second variable.

Let us moreover assume that $X$ is projective.
A $v$-adic metric on the trivial line bundle is characterized
by the norm of its section~$1$, which is 
a nonvanishing continuous function~$\rho_v$ on~$X(F_v)$.
Consequently, an adelic metric~$\rho$ on~$\mathscr O_X$
is given by a family~$(\rho_v)$ of such functions
which are almost all equal to~$1$. (The projectivity
assumption on~$X$ allows to write $X$ as a finite union of open subsets~$U$
such that $X(F_v)$ is the union of the sets $U(\mathfrak o_v)$.)
One has
\[ H_\rho(x) = \prod_{v\in\Val(F)} \rho_v(x)^{-1} \]
(finite product).
If $\rho_v$ is identically equal to~$1$, we set $c_v=1$.
Otherwise, by continuity and compactness of $X(F_v)$, there is a positive
real number~$c_v$ such that $c_v^{-1}<\rho_v(x)<c_v$ for any $x\in X(F_v)$.
Let $c$ be the product of all~$c_v$ (finite product);
it follows that 
$ c^{-1}< H_\rho(x)< c$ for any $x\in X(F)$.

As a consequence, if $L$ and~$L'$ are two adelically metrized line bundles 
with  the same underlying line bundle, the quotient of the
exponential height functions $H_L$ and~$H_{L'}$ is bounded
above and below (by a positive real number). This property
is best expressed using the logarithmic heights $h_L=\log H_L$
and $h_{L'}=\log H_{L'}$: it asserts that the difference
$h_L-h_{L'}$ is bounded.

\subsubsection{Functoriality}
If $f\colon X\ra Y$ is a morphism of algebraic varieties
and $L$ is a adelically metrized line bundle  on~$Y$,
the line bundle $f^*L$ on~$X$ has a natural adelic metric:
indeed, the fibre of $f^*L$ at a point~$x$ is the line $L(f(x))$,
for which we use its given norm. This construction is compatible
with isometry classes and defines a morphism of Abelian groups
$f^*\colon\hPic(Y)\ra\hPic(X)$.
At the level of heights, it implies the equality
\[ H_L(f(x))=H_{f^*L}(x) \]
for any $x\in X(F)$.

\subsubsection{Finiteness property}
We have explained that for any real number~$B$,
$\P^n(\Q)$ has only finitely many points of height smaller than~$B$.
As was first observed by~\cite{northcott1950},
this property extends to our more general setting:
\emph{if $L$ is an adelically metrized line bundle whose
underlying line bundle is ample, then the set of points $x\in X(F)$
such that $H_L(x)\leq B$ is finite.}
 
This can be proved geometrically on the basis of the result for
the projective space over~$\Q$. For simplicity of notation,
we switch to logarithmic heights.
Since $L$ is ample, it is known
that there exists an integer~$d\geq 1$ such that $L^{\otimes d}$
is very ample,
hence of the form $\phi^*\mathscr O(1)$ for some
embedding $\phi$ of~$X$ into a projective space~$\P^n$.
The functoriality property of heights together with the boundedness
of the height when the underlying metrized line bundle is trivial,
imply that there is a positive real number~$c$ such that
\[ h_L(x)=\frac 1d h_{L^{\otimes d}}(x) \geq \frac1d h(\phi(x))+c. \]
Consequently, it suffices to show the finiteness result on~$\P^n$.
Moreover, the explicit formula given for the height on~$\P^n$ shows that
\[ h([x_0:\cdots:x_n)) \geq \sup_{i,j} h([x_i:x_j]), \]
where the supremum runs over all couples~$(i,j)$ such that $x_i$
and~$x_j$ are not both equal to~$0$.
This reduces us to proving the finitess assertion on~$\P^1$.

We also may assume that $F$ is a Galois extension of~$\Q$ and
let $G=\{\sigma_1,\dots,\sigma_d\}$ be its Galois group.
Let $x\in F$ and let us consider the polynomial 
\[ P_x=\prod_{\sigma\in G}(T-\sigma(x))=T^d+a_1T^{d-1}+\dots+a_d;\]
it has coefficients in~$\Q$.
For $v\in\Val(F)$ and $\sigma\in G$, the function
$t\mapsto \abs{\sigma(t)}_v$ is a normalized absolute value on~$F$
which is associated to the $p$-adic  absolute value if
so is~$v$. Consequently,
\[ \max(1,\abs{\sigma_1(x)}_v,\dots,\abs{\sigma_d(x)}_v)
 \leq \prod_{w|p} \max(1,\abs x_w),  \]
where the notation~$w\mid p$ 
means that the product ranges over all normalized absolute values~$w$ on~$F$
which are equivalent to the $p$-adic absolute value.
It follows that $h([1:\sigma_1(x):\dots:{\sigma_d(x)})\leq d h([1:x])$.

Moreover, the map which associates to a point $[u_0:\dots:u_d]\in\P^d$
the coefficients $[v_0:\dots:v_d]$ of the polynomial $\prod_{i=1}^d (u_0T-u_i)$
has degree~$d$. This implies that
\[ h([v_0:\dots:v_d]) = d h([u_0:\dots:u_d])+\mathrm O(1). \]
If $h(x)\leq T$, this implies that the point $\mathbf a=[1:a_1:\dots:a_d]$
of~$\P^d(\Q)$ is of height $h(\mathbf a)\leq d^2 T$.
By the finiteness property over~$\Q$, there are only finitely
many possible polynomials~$P_x$, hence finitely many~$x$
since $x$ is one of the~$d$ roots of~$P_x$.

\subsubsection{Rational points \emph{vs} algebraic points}
To keep the exposition at the simplest level, we have only
considered the height of $F$-rational points. However,
there is a suitable notion of adelic metrics
which uses not only the $v$-adic points $X(F_v)$,
but all points of~$X$ over the algebraic closure $\bar{F_v}$ of~$F_v$,
or even its completion~$\C_v$. Together with the obvious  extension 
to metrized line bundles of the relation between $H_F(\mathbf x)$ and~$H_{F'}(\mathbf x)$ when $F'$ is a finite extension of~$F$ and $\mathbf x\in\P^n(F)$
this allows to define a (exponential) height function $H_L$ on the whole of $X(\bar F)$, as well as its logarithmic counterpart $h_L=\log H_L$.
The preceding properties extend to this more general setting,
the only nonobvious assertion being the boundedness of the height
when the underlying line bundle is trivial.

\section{Manin's problem}
\subsection{Counting functions and zeta functions}
\subsubsection{The counting problem}
Let $X$ be a variety over a number field~$F$
and let $L$ be an ample line bundle with an adelic metric.
We have seen that for any real number~$B$, the set of points
$x\in X(F)$ such that $H_L(x)\leq B$ is a finite set.
Let $N_X(L;B)$ be its cardinality.
We are interested in the asymptotic behaviour of $N_X(L;B)$,
when $B$ grows to infinity. We are also interested in understanding
the dependence on~$L$.

\subsubsection{Introducing a generating series}
As is common practice in all counting problems, \emph{e.g.}, in
combinatorics, one introduces a generating series.  For the present
situation, this is a \emph{Dirichlet series}, called the
\emph{height zeta function}, and defined by
\[ Z_X(L;s) = \sum_{x\in X(F)} H_L(x)^{-s}, \]
for any complex number $s$ such that the series converges absolutely.
In principle, and we will eventually do so,
one can omit the ampleness condition in that definition,
but then the convergence in some half-plane is not assured
(and might actually fail).

\subsubsection{Abscissa of convergence}
The abscissa of convergence $\beta_X(L)$ is the infimum
of all real numbers~$a$ such that $Z_X(L;s)$ converges absolutely
for $\Re(s)>a$.
Since the height zeta function is a Dirichlet series with positive
coefficients, a theorem of Landau implies that $\beta_X(L)$
is also the infimum of all real numbers~$a$ such that $Z_X(L;a)$
converges.

Let us assume again that $L$ is ample.
Then, the proof of \textsc{Northcott}'s theorem shows that $N_X(L;B)$
grows at most polynomially. It follows that $Z_X(L;s)$
converges for $\Re(s)$ large enough. For instance,
let us assume that $N_X(L;B)\ll B^a$. Let us fix a real number $B>1$.
There are $\ll B^{na}$ points $x\in X(F)$
such that $H_L(x)\leq B^n$ and the part of the series
given by points of heights between~$B^{n-1}$ and~$B^{n}$
is bounded by $B^{na} B^{-(n-1)s}=B^{n(a-s)+s}$. 
By comparison with the geometric
series, we see that $Z_X(L;s)$ converges for   $\Re(s)>a$.

This shows that the function~$L\mapsto \beta_X(L)$ is well defined
for adelically metrized line bundles whose underlying 
line bundle is ample.
Moreover, it does not depend on the choice of an adelic metric,
so it really comes from a function on the set of ample
line bundles in~$\Pic(X)$. 
The formula $H_{L^{\otimes d}}(x)=H_L(x)^d$ implies
that $\beta_X(L^{\otimes d})=d\beta_X(L)$ for any positive integer~$d$.
In other words, $\beta_X$ is homogeneous of degree~$1$.
One can also prove that it uniquely extends to a continuous
homogeneous function on~$\Lambda_\ample$
(see \cite{batyrev-m90}).

\subsubsection{Analytic properties of $Z_X(L;\cdot)$ \emph{vs.}
asymptotic expansions of $N_X(L;\cdot)$}
Despite this simple definition, very little is known about $\beta_X$,
let alone about the height zeta function itself.
However, there are numerous examples where $Z_X(L;\cdot)$
has a meromorphic continuation to some half-plane,
with a unique pole of largest real value, say $\alpha_X(L)$,
of order~$t_X(L)$. By \textsc{Ikehara}'s Tauberian theorem,
this implies an asymptotic expansion of the form
\[ N_X(L;B) \sim c B^{\alpha_X(L)} (\log B)^{t_X(L)-1}. \]
In some cases, one can even establish terms of lower order,
or prove explicit error terms.

However, these zeta functions are not as well behaved as
the ones traditionally studied in algebraic number theory.
Although some of them have a meromorphic extension to the whole complex
plane, like those of flag varieties (see \S\ref{sec.eisenstein} below,
see also~\cite{essouabri2005}),
it happens quite often, for example, that they have a natural boundary ;
this is already the case for some toric surfaces,
see~\cite[Example~3.5.4]{batyrev-t95b}.
Subtler analytic properties
of the height zeta functions beyond the largest pole
is the subject of some recent investigations,
see, \emph{e.g.}, \cite{breteche-swd2007}.

\subsection{The largest pole}

\subsubsection{The influence of the canonical line bundle}
We now assume that $X$ is smooth.
One of the insights of mathematicians in the \textsc{xx}th century
(notably \textsc{Mordell}, \textsc{Lang} and~\textsc{Vojta})
was that the potential \emph{density}
of rational points is strongly related by the negativity of the
canonical line bundle~$\omega_X$ with respect to the ample
cone.

This is quite explicit in the case of curves. Namely, if $X$
is a curve of genus~$g$, three cases are possible:
\begin{enumerate} 
\item \textbf{genus $g=0$, $\omega_X^{-1}$ ample.}
Then $X$ is a \emph{conic}. Two subcases are possible:
either $X(F)$ is empty---$X$ has no rational point---
or $X$ is isomorphic to the projective line. Moreover,
the \textsc{Hasse} principle allows to decide quite effectively
in which case we are.
\item  \textbf{genus $g=1$, $\omega_X$ trivial.}
If $X$ has a rational point, then $X$ is an \emph{elliptic curve},
endowing $X(F)$ with a structure of Abelian group. Moreover,
the theorem of \textsc{Mordell-Weil} asserts that
$X(F)$ is of finite type.
\item \textbf{genus $g\geq 2$, $\omega_X$ ample.}
By \textsc{Mordell}'s conjecture, first proved  by~\cite{faltings83},
$X(F)$ is a finite set.
\end{enumerate}

The counting function distinguishes very clearly
the first two cases. When $X$ is the projective line,
$N_X(B)$ grows like a polynomial in~$B$ (depending on how the height
is defined), and for $X$ an elliptic curve, $N_X(B)\approx (\log B)^{r/2}$,
where $r$ is the rank of the Abelian group~$X(F)$.

\subsubsection{Increasing the base field}
The general conjectures of~\textsc{Lang} lead us to expect
that if $\omega_X$ is ample, then $X(F)$ should not be
dense in~$X$, but that this could be expected in
the opposite case where $\omega_X^{-1}$ is ample.

However, as the case of conics already indicates, geometric invariants
like the ampleness of~$\omega_X^{-1}$ cannot suffice to decide
on the existence of rational points. For example, the conic~$C$ over~$\Q$
given by the equation $x^2+y^2+z^2=0$ in the projective plane~$\P^2$
has no rational points (a sum of three squares of rational
integers cannot be zero,
unless they are all zero). However, as soon as the ground field~$F$
possesses a square root of~$-1$, then~$C$ admits a rational
point~$P=[1:\sqrt{-1}:0]$ and the usual process of intersecting
with the conic~$C$ a variable line through that point~$P$ gives
us a parametrization of~$C(F)$.

Smooth projective varieties with~$\omega_X^{-1}$ ample
are called \emph{Fano varieties}. In general,
for such a variety over  a number field~$F$, it is expected 
that \emph{there exists a finite extension~$F'$ of~$F$
such that $X(F')$ is dense in~$X$ for the Zariski topology.}
This seems to be a very difficult question to settle in general,
unless $X$ has particular properties, like being rational or unirational,
in which case a dense set of points in~$X(F')$ can be parametrized
by the points of a projective space~$\P^n(F')$ (where $n=\dim X$).

\subsubsection{Eliminating subvarieties}
This conjecture considers the ``density aspect'' of the rational points.
Concerning the counting function, the situation is even more complicated
since nothing forbids that most of the rational points of small height
are contained in subvarieties.

Let $Y$ be a subscheme of~$X$. One says that $Y$ is \emph{strongly accumulating}
if the fraction  $N_Y(L,B)/N_X(L,B)$ tends to~$1$ when $B\ra\infty$,
and one says that $Y$ is \emph{weakly accumulating} if the inferior
limit of this fraction is positive.

As a consequence, the behaviour of the counting function
can only be expected to reflect the global geometry if one
doesn't count points in accumulating subvarieties.
This leads to the notation $N_U(L;B)$
and $Z_U(L;s)$, where $U$ is any Zariski open subset in~$X$.

To get examples of such behaviour, it suffices to blow-up
a variety at one (smooth) rational point~$P$.
It replaces the point~$P$ by a projective space~$E$ of dimension~$n-1$
if $X$ has dimension~$n$ which is likely to possess $\approx B^a$
points of height~$\leq B$, while the rest could be smaller.
For example, if $X$ is an Abelian variety (generalization of
elliptic curves in higher dimension), $N_X(L;B)$ only grows
like a power of~$\log B$.

\subsubsection{Definition of a geometric invariant}
In order to predict on a geometric basis the abscissa
of convergence of $Z_U(L;s)$, one needs to introduce
a function on~$\Pic(X)_\R$ which is homogeneous of degree~$1$
and continuous. It is important to know at that point
that for Fano varieties, this vector space~$\Pic(X)_\R$
is finite-dimensional (and identifies with the so-called
Néron--Severi group).
Indeed, since $\omega_X^{-1}$ is ample
(this is the very definition of a Fano variety),
\textsc{Serre}'s duality and \textsc{Kodaira}'s vanishing theorem
imply that $\mathrm H^1(X,\mathscr O_X)=0$, so that
the Albanese variety of~$X$ is trivial.

After considering many examples (\cite{schanuel79}, \cite{serre1997},
\cite{franke-m-t89},...), and with detailed investigations of surfaces at hand,
\cite{batyrev-m90} concluded that
the relevant part of $\Pic(X)_\R$ was not the ample cone
(as the case of the projective space naïvely suggests),
but the \emph{effective cone}~$\Lambda_\eff$.

In effect, they define for any line bundle~$L$
belonging to the interior of the effective cone
a real number~$\alpha_X(L)$ which is the 
least real number~$a$ such that $\omega_X\otimes L^{\otimes a}$
belongs to the effective cone in~$\Pic(X)_\R$.

\subsubsection{The conjectures of Batyrev and Manin}
The consideration of the effective cone fits nicely with
the elimination of accumulating  subvarieties. Indeed, 
let us assume that the line bundle $L$ belongs  to the \emph{interior}
of the effective cone; then there exists a positive integer~$d$,
an ample line bundle~$M$ and an effective line bundle~$E$
on~$X$ such that $L^{\otimes d}\simeq M\otimes E$.
(These line bundles are also called \emph{big}.)
For any nonzero section~$s_E$ of~$E$, $H_M(x)$ is bounded from
below where $s_E$ does not vanish, so that $H_L(x)\gg H_M(x)^{1/d}$
on the open set $U=X\setminus\{s_E=0\}$. 
It follows that $N_U(L;B)$ grows at most polynomially,
and that the abscissa of convergence $\beta_U(L)$ of $Z_U(L;s)$ is finite.

\cite{batyrev-m90} present three conjectures of increasing precision
concerning the behaviour of the counting function.
Let $X$ be a projective smooth variety over a number field~$F$,
let $L$ be a line bundle on~$X$ which belongs to the interior
of the effective cone. The following assertions are conjectured:
\begin{enumerate}
\item For any $\eps>0$, there exists a dense open subset~$U\subset X$ such
that $\beta_U(L)\leq \alpha_X(L)+\eps$
(\emph{linear growth conjecture}).
\item If $X$ is a Fano variety, then for any large enough finite
extension~$F'$ of~$F$, and any small enough dense open set~$U\subset X$,
one has $\beta_{U,F'}(L)=\alpha(L)$.
\item Same assertion, only assuming that $\omega_X$
does not belong to the effective cone.
\end{enumerate}

\subsection{The refined asymptotic expansion}
\subsubsection{Logarithmic powers}
In the cases evoked above, the height zeta function appears to possess
a meromorphic continuation to the left of the line $\Re(s)=\alpha_X(L)$,
with a pole of some order~$t$ at $s=\alpha_X(L)$. By Tauberian theory,
this implies a more precise asymptotic expansion for the
counting function, namely 
$N_U(L;B)\approx B^{\alpha_X(L)} (\log B)^{t-1}$.
A stronger form of the conjecture
of \textsc{Batyrev--Manin} also predicts  the order of this pole,
at least when $L=\omega_X^{-1}$. In that case, \textsc{Batyrev}
and \textsc{Manin} conjecture that $t$
is the dimension of the real vector space~$\Pic(X)_\R$.
(Recall that the Picard group of a Fano variety has finite rank.)

\subsubsection{(In)compatibilities}
Some basic checks concerning that conjecture had been made
by~\cite{franke-m-t89}, for instance its compatibility with products
of varieties.
However, \cite{batyrev-t96b} observed that \emph{the conjecture
is not compatible with families} and produced a \emph{counterexample}.

They consider the subscheme~$V$ of~$\P^3\times\P^3$ 
defined by
the equation $x_0y_0^3+\cdots+x_3y_3^3=0$
(where $[x_0:\dots:x_3]$ are the homogeneous coordinates
on the first factor, and $[y_0:\dots:y_3]$ are those on the second).
For fixed~$\mathbf x\in\P^3$, the fibre~$V_{\mathbf x}$
consisting of $\mathbf y\in\P^3$ such that $(\mathbf x,\mathbf y)\in V$
is a diagonal cubic surface; in other words, $V$ is the total
space of the family of diagonal cubic surfaces.

Let $p_1$ and $p_2$ be the two projections from~$\P^3\times\P^3$ to~$\P^3$.
For any couple of integers~$(a,b)$, let $\mathscr O(a,b)$ be
the line bundle $p_1^*\mathscr O_{\P^3}(a)\otimes p_2^*\mathscr O_{\P^3}(b)$
on~$\P^3\times\P^3$.
The so-defined map from~$\Z^2$ to~$\Pic(\P^3\times\P^3)$ is an isomorphism
of Abelian groups. Moreover, since $V$ is an hypersurface in~$\P^3\times\P^3$,
the Lefschetz theorem implies that the restriction map induces an isomorphism
from $\Pic(\P^3\times\P^3)$ to~$\Pic(V)$.
The anticanonical line bundle~$\omega_V^{-1}$ of~$V$ 
corresponds to~$\mathscr O(3,1)$ (anticanonical of~$\P^3\times \P^3$
minus the class of the equation, that is $(4,4)-(1,3)$)
and  the conjecture of \textsc{Batyrev--Manin} predicts that there
should be $\approx B(\log B)$ points of $\mathscr O(3,1)$-height~$\leq B$
in~$V(F)$.

Let us fix some~$\mathbf x\in\P^3$ such that $V_{\mathbf x}$ is non-singular.
This is a cubic surface, and its anticanonical line bundle
identifies with the restriction of~$\mathscr O_{\P^3}(1)$.
Thus, if $t_{\mathbf x}=\dim_\R\Pic(V_{\mathbf x})_\R$, 
the conjecture predicts that there should be
$\approx B(\log B)^{t_{\mathbf x}-1}$ points of $\mathscr O(1)$-height~$\leq B$
in~$V_{\mathbf x}(F)$.

As a consequence, if the conjecture holds for~$V$
and if  $t_{\mathbf x}>2$, then the  fibre~$V_{\mathbf x}$ would have
more points than what the total space seem to have.
This doesn't  quite contradict the conjecture however, because it 
explicitly allows to remove a finite union of subvarieties.

However, the rank of the Picard group can exhibit disordered behaviour in 
families; for example, it may not be semi-continuous,
and jump on a infinite union  of subvarieties.
This happens here, since $t_{\mathbf x}=7$ 
when $F$ contains the cubic roots of unity
and all the homogeneous coordinates of~$\mathbf x$ are cubes.
The truth of the conjecture for~$V$ therefore requires to
omit all such fibres~$V_{\mathbf x}$. But they form an infinite
union of disjoint subvarieties, a kind of accumulating subset
which is not predicted by the conjecture of \textsc{Batyrev--Manin}.

In particular, this conjecture is false, either for~$V$,
or for most cubic surfaces. 
By geometric considerations, and using their previous
work on toric varieties, \cite{batyrev-t96b} could in fact
conclude that the conjecture does not hold for~$V$.

\subsubsection{Peyre's refinement of the conjecture}
One owes to~\cite{peyre95} a precise refinement of the previous
conjecture, as well as the verification of the refined conjecture
in many important cases.
Indeed, all known positive examples feature an asymptotic expansion
of the form $N_U(\omega_X^{-1};B)\sim c B(\log B)^t$,
for some positive real number~$c$, which 
in turn is the product of four factors:
\begin{enumerate}
\item the volume of a suitable subspace $X(\AD_F)^{\BM}$
of the adelic space $X(\AD_F)$
(where the so-called Brauer-Manin obstruction to rational 
points vanishes) with respect to the \emph{Tamagawa measure}
first introduced in that context by~\cite{peyre95}.
\item the cardinality of the finite Galois cohomology group 
$\mathrm H^1 (\Pic(X_{\bar F}))$; 
\item  a rational number related to the location of the anticanonical line 
bundle in the effective cone $\Lambda_\eff$ of~$\Pic(X)_\R$;
\item the rational number $1/t !$. 
\end{enumerate}

Let us give a short account of \textsc{Peyre}'s construction of a measure.
The idea is to use the given adelic metric on $\omega_X^{-1}$
to construct a measure~$\tau_{X,v}$ on the local spaces $X(F_v)$
for all places~$v$ of the number field~$F$, 
and then to consider a suitable (renormalized) product of these
measures.

\subsubsection{Local measures}
\label{subsubsec.local-measures.I}
Ideally, a measure on an $n$-dimensional analytic manifold is given in local
coordinates by an expression of the form 
$\mathrm d\mu(x)=f(x) \,\mathrm dx_1\dots \mathrm dx_n$,
where $f$ is a positive function. When one considers another
system of local coordinates, the expression of the measure
changes, and its modification is dictated by the change-of-variables
formula for multiple integrals: the appearance of the absolute
value of the Jacobian means that $\mathrm d\mu(x)$ is modified
as if it were ``the absolute value of a differential form''.
One can deduce from this observation that if $\alpha$ is a local
(\emph{i.e.}, in a chart) differential $n$-form, written
$\alpha =f(x) \mathrm dx_1\wedge\dots\wedge\mathrm dx_n$
in coordinates, then the measure  
\[ \frac{\abs\alpha_v}{\norm\alpha_v} = \frac{\abs{f(x)}_v}{\norm \alpha_v}
    \,\mathrm dx_1\dots \mathrm dx_n \]
is well-defined, independently of the choice of~$\alpha$.
One can therefore glue all these local measures
and obtain a measure~$\tau_{X,v}$ on~$X(F_v)$.

\subsubsection{Convergence of an infinite product}
Since $X$ is proper,
the adelic space $X(\AD_F)$ is the (infinite) product of all $X(F_v)$,
for $v\in\Val(F)$, endowed with the product topology.
To define a measure on~$X(\AD_F)$ one would like to consider
the (infinite) product of the measures~$\tau_{X,v}$.
However, this gives a finite measure if and only if the product
\[ \prod_{\text{$v$ finite}} \tau_{X,v}(X(F_v)) \]
converges absolutely.
This convergence never holds, 
and one therefore needs to introduce \emph{convergence factors}
to define a measure on $X(\AD_F)$.

By a  formula of~\cite{weil82},
the equality $\tau_{X,v}(X(F_v))=q_v^{-\dim X}\Card (X(k_v))$
holds for almost all finite places~$v$,
where $k_v$ is the residue field of~$F$ at~$v$, $q_v$
is its cardinality, and $X(k_v)$ is the set of solutions in~$k_v$
of a fixed system of equations with coefficients in the
ring of integers of~$F$ which defines~$X$.
By \textsc{Weil}'s conjecture, established by~\cite{deligne74}, plus various
cohomological computations, this implies that
\[ \tau_{X,v}(X(F_v)) = 1+ \frac 1{q_v} \Tr(\Frob_v|\Pic(X_{\bar F})_\R)
 + \mathrm O(q_v^{-3/2}), \]
where $\Frob_v$ is a ``geometric Frobenius element'' 
of the Galois group of~$\bar F$  at the place~$v$.

Let us consider the Artin L-function of the 
Galois-module~$P=\Pic(X_{\bar F})_\R$:
it is defined as the infinite product
\[ \mathrm L(s,P)=\prod_{\text{$v$ finite}} \mathrm L_v(s,P),
\qquad \mathrm L_v(s,P) = \det(1-q_v^{-s}\Frob_v|P)^{-1}. \]
The product converges absolutely for $\Re(s)>1$,  defines a holomorphic 
function in that domain. Moreover, $\mathrm L(s,P)$ has a meromorphic
continuation to~$\C$ with a pole of order~$t=\dim\Pic(X)_\R$
at $s=1$. Let 
\[ \mathrm L^*(1,P)=\lim_{s\ra 1} \mathrm L(s,P) (s-1)^{-t}; \]
it is a positive real number.

Comparing the asymptotic expansions for $\mathrm L_v(1,P)$
and $\tau_{X,v}(X(F_v))$, one can conclude that the following
infinite product
\[ \tau_X = \mathrm L^*(1,P) \prod_v \left( \mathrm L_v(1,P)^{-1}\tau_{X,v} \right)
\]
is absolutely convergent
and defines a Radon measure on the space $X(\AD_F)$.

\subsubsection{Equidistribution}
One important aspect of the language of metrized line bundles is that
it suggests explicitly to look at what
happens when one changes the adelic metric on the 
canonical line bundle.

Assume for example that there is a number field~$F$ and an open subset~$U$
such that \textsc{Manin}'s conjecture (with precised constant as above) holds
for \emph{any} adelic metric on~$\omega_X$.
Then, \cite{peyre95} shows that rational points in~$U(F)$ of heights~$\leq B$ 
\emph{equidistribute} towards the probability measure on~$X(\AD_F)^\BM$
proportional to~$\tau_X$. 
Namely, for any smooth function~$\Phi$ on $X(\AD_F)$,
\[ \frac1{N_U(\omega_X^{-1};B)} 
 \sum_{\substack{x\in U(F)\\ H_{\omega_X^{-1}}(x)\leq B}} \Phi(x)
   \longrightarrow \frac1{\tau_X(X(\AD_F)^\BM)} \int_{X(\AD_F)^\BM} \Phi(x) \,\mathrm d\tau_X(x).\]

As we shall see in the next section, this strengthening
by \textsc{Peyre} of the conjecture of \textsc{Batyrev--Manin}
is true in many remarkable cases, with quite nontrivial proofs.

\section{Methods and results}

\subsection{Explicit counting}

\subsubsection{Projective space}
Let $F$ be a number field and let $L$ be the line bundle $\mathscr O(1)$
on the projective space~$\P^n$, with an adelic metric.
\cite{schanuel79} has given the following asymptotic expansion
for the counting function on~$\P^n$:
\[ N_{\P^n}(L;B)  \sim \frac{\Res_{s=1}\zeta_F(s)}{\zeta_F(n+1)} 
\left( \frac{2^{r_1}(2\pi)^{r_2}}{\sqrt{\abs{D_F}}} \right)^n (n+1)^{r_1+r_2-1}
B^{n+1}. \]
In this formula, $\zeta_F$ is the Dedekind zeta function of the
field~$F$, $r_1$ and~$2r_2$ are the number of real embeddings
and complex embeddings of~$F$, and $D_F$ is its discriminant.

In fact, it has been shown later by many authors, 
see \emph{e.g.} \cite{franke-m-t89},
that the height zeta function $Z_{\P^n}(L;s)$ is holomorphic
for $\Re(s)>n+1$, has a meromorphic continuation to the whole
complex plane~$\C$, with a pole of order~$1$ at the point~$s=n+1$
and no other pole on the line $\{\Re(s)=n+1\}$.

To prove this estimate, one can sort points $[x_0:\dots:x_n]$
in~$\P^n(F)$ according to the class of the fractional ideal 
generated by~$(x_0,\dots,x_n)$. Constructing fundamental domains
for the action of units,
the enumeration of such sets can be reduced to the counting
of lattice points in homothetic sets in~$\R^{(n+1)(r_1+2r_2)}$
whose boundary has smaller dimension. The Möbius inversion formula
is finally used to take care  of the coprimality condition.

\subsubsection{Hirzebruch surfaces and other blow-ups}
There are many other cases where this explicit method
works. Let us mention only some of them:
\begin{itemize}
\item ruled surfaces over the projective line (Hirzebruch surfaces);
\item Grassmann varieties (\cite{thunder92}) ;
\item Del Pezzo surfaces of degree~$5$, 
given as blow-ups of the plane in 4 points in general position
(\cite{breteche2002});
\item some Chow varieties of~$\P^n$, see~\cite{schmidt1993b}.
\end{itemize}

This method suffers from an obvious drawback: 1) it is hard to take
advantage of subtle geometric properties of the situation studied;
2) it requires to \emph{know} already much about the rational points,
\emph{e.g.}, to be able to parametrize them.
However, this is so far the only way of dealing with
the (eventually singular) Del Pezzo surfaces which are not 
equivariant compactifications of the torus $\gm^2$
or of the additive group~$\ga^2$.

To deal with the parametrization of the rational points,
an essential tool is the \emph{universal torsor},
invented by \cite{colliot-thelene-sansuc1987},
which is a quasi-projective variety lying \emph{over}
the given Del Pezzo surface.
The computation of a universal torsor
is essentially equivalent to that of a multi-graded
ring defined in~\cite{cox1995}, the so-called \emph{Cox ring}
of the variety. 
In the case of the projective space~$\P^n$, the torsor
is $\A^{n+1}\setminus\{0\}$, mapping to~$\P^n$
via $(x_0,\ldots,x_n)\mapsto [x_0:\dots:x_n]$,
and the Cox ring is the graded polynomial ring
in~$n+1$ variables.

In the context of \textsc{Manin}'s conjecture,
universal torsors were first introduced by~\cite{salberger98} 
for toric varieties, in which case they are isomorphic
to an open subset of an affine space.
The parametrization of rational points which  they allow to write
proved to be very useful for the study of points
of bounded height, at least in contexts where other techniques, like
harmonic analysis, do not apply.
Moreover, their computation is also related to the precise
understanding of the numerically effective cone in
the Picard group. 

\cite{derenthal2006} computes these universal torsors 
for smooth Del Pezzo surfaces
(that is, blow-ups of the projective plane in at most~$8$
points in general position).
They were also computed for many singular Del Pezzo surfaces,
see, \emph{e.g.}, the article of~\cite{derenthal-tschinkel2007}
for an example leading to a case of \textsc{Manin}'s conjecture.
In that paper, the interested reader shall also find
a  list of known cases, mostly due to 
\textsc{Browning}, de la \textsc{Bretèche} and~\textsc{Derenthal},
with references to the articles; 
see also the more recent preprint of~\cite{derenthal2008}.
Anyway, we apologize not to develop this important 
part of the subject here, which would lead us to far of
your topic.

\subsection{The circle method of Hardy--Littlewood}
The \emph{circle method} was initially devised  to tackle
\textsc{Waring}'s problem, namely the decomposition of an integer
into sums of powers. More generally, it is well suited to the study
of diophantine equations in ``many variables''.

Concerning our counting problem,
it indeed allows to  establish the conjecture of 
\textsc{Batyrev--Manin--Peyre}
for smooth complete intersections of codimension~$m$ in~$\P^n$,
that is, subschemes~$X$ of~$\P^n$ defined by the vanishing of $m$ homogeneous 
polynomials $f_1,\dots,f_m$ such that at every point~$x\in X$,
the Jacobian matrix of the~$f_i$ has rank~$m$.
Let $d_1,\dots,d_m$ denote the degrees of the polynomials~$f_1,\dots,f_m$;
then the canonical line bundle of~$X$ is the restriction to~$X$
of the line bundle~$\mathscr O(-n-1+\sum_{i=1}^m d_i)$. Consequently,
$X$ is Fano if and only if $d_1+\dots+d_m\leq n$.
Therefore, specific examples of Fano complete intersections are lines or conics
in~$\P^2$, planes, quadrics or cubic surfaces in~$\P^3$, etc.

The circle method is restricted to cases where the number of variables~$n$
is very large in comparison to the degrees~$d_1,\dots,d_m$.
For instance,  \cite{birch62} establishes the desired result for an hypersurface
of degree~$d$ in~$\P^n$ provided $n\geq 2^d(d-1)$.
However, when it applies, it gives very strong results,
including a form of the Hasse principle and equidistribution,
see~\cite{peyre95}.

\subsection{Homogeneous spaces}\label{sec.homogeneous}

Homogeneous spaces, namely quotients of an algebraic group
by a subgroup, form a large and interesting class of varieties.
Here is a list of interesting examples:
\begin{itemize}
\item
Algebraic groups themselves, like $\SL_n$ viewed as the subset
of~$\P^{n^2}$ defined by the equation  $z^n=\det(A)$,
where a point of~$\P^{n^2}$ is a non-zero
pair~$(z,A)\in F\times \operatorname{Mat}_n(F)$
modulo homotheties.
\item
The set $V_P$ of matrices with fixed characteristic polynomial~$P$.
Indeed, when the polynomial~$P$ 
has distinct roots, any two such matrices are conjugate,
hence $V_P$ is an homogeneous space of~$\GL_n$.
\item
Non singular quadrics;
indeed, if $Q$ is a non-degenerate quadratic form in~$n$ variables,
the theorem of Witt implies
that the orthogonal group~$\mathrm O(Q)$ acts transitively
on the hypersurface of~$\P^{n-1}$ defined by~$Q$.
\item
Grassmann varieties: if $G$ is the Grassmann variety of $d$-planes~$W$
in the $n$-dimensional affine space~$\A^n$, it has a transitive
action of the group~$\GL(n)$, given by $W\mapsto gW$;
the stabilizer of the fixed $d$-plane~$\A^d\times\{0\}$
is isomorphic the subgroup of block-triangular matrices (see below).
\end{itemize}

Various techniques allow to tackle the question
of rational points of bounded height on such varieties.
These methods allow to derive asymptotic expansions
for the counting function by comparing it to the volume of a corresponding
set in a real or adelic space, whose asymptotic growth
has to be established independently. They are also amenable to proving
equidistribution results, as well as establishing
analyting properties of the height zeta function. 

\subsubsection{Point counting \emph{vs.} balls}
Three important  papers appeared in the 90's
concerning the close problem of integral points of bounded height
and still have a great influence.
Methods of ergodic theory apply usually when one considers
\emph{lattices} of real or adelic semi-simple groups
(that is, discrete subgroups of finite covolume).
Keywords are the equidistribution theorem of M.~\textsc{Ratner}
for unipotent flows (used in the closely related context
of integral points by~\cite{eskin-m-s1996}) or,
as in~\cite{eskin-mcmullen1993,gorodnik-m-o2009},
mixing properties of the
dynamical system, implied by decay of matrix coefficients
of unitary representations (theorem of~\textsc{Howe--Moore}
and its extensions).

\subsubsection{Eisenstein series}\label{sec.eisenstein}
\emph{Flag varieties} are classical
generalizations of the projective space:
they parametrize subvector spaces~$W$ of a fixed vector space~$V$
(Grassmann variety), of, more generally, increasing
families $W_1\subset\cdots\subset W_m$ (``flags'') of subvector spaces.
Let us consider for simplicity the case of a fixed vector space~$V$
of dimension~$n$
and the Grassmann variety $\Gr_d^n$ of subspaces of fixed dimension~$d$.
The projective space~$\P(V)$ is recovered by considering
the cases $d=1$ (parametrizing lines) and $d=n-1$
(parametrizing hyperplanes).  
Let us fix a specific subspace~$W_0$ of dimension~$d$ in~$V$; 
considering a basis of~$W_0$ and extending it to a basis of~$V$,
it is easy to observe that any other subspace of dimension~$d$
is of the form $gW_0$, for some automorphism $g\in\GL(V)$;
moreover, $gW_0=g'W_0$ if and only if $g^{-1}g'W_0=W_0$,
that is, if $g^{-1}g'$ belongs to the \emph{stabilizer}~$P$ of~$W_0$
for the action of~$\GL(V)$. In the basis that we fixed, 
$P$ is the set of invertible upper-triangular block matrices having
only zeroes in the lower rectangle $[d+1,\dots,n]\times[1,\dots,d]$.

More generally, in the language of algebraic groups, generalized
flag varieties are quotients of reductive groups  by parabolic subgroups.
The understanding of their height zeta functions is a consequence
of the observation, due to J.~\textsc{Franke}, that
these functions were  studied extensively in the theory of automorphic forms,
being nothing else than generalized \emph{Eisenstein series}.
Modulo some (not so obvious) translation,
the results of  R.~P.~\textsc{Langlands} on these Eisenstein 
series readily imply the conjectures of \textsc{Batyrev--Manin}
and \textsc{Peyre} (including equidistribution) for the flag varieties,
see~\cite{franke-m-t89} and~\cite{peyre95}.

\subsubsection{Spectral decomposition of the height zeta function}
Beyond Eisenstein's series and \textsc{Langlands}'s results,
there are many interesting cases  where the
spectral decomposition of the height zeta function
can be computed, leading to a proof of \textsc{Manin}'s conjecture.
In fact, this  has been done for equivariant compactifications of 
many families of algebraic groups.

To give an idea of the principle of the proof, let us say a word
on the analogous problem of counting lattice points in 
large balls. Let $B$ be a bounded measurable subset in~$\R^n$
and let $\chi_B$ be its characteristic function.
To understand the cardinality~$N(B)$ of $B\cap \Z^n$,
\emph{i.e.}, the sum $\sum_{x\in\Z^n}\chi_B(x)$, one can
try to use the Poisson summation formula which would relate
this sum to a sum over $\Z^n$ of the Fourier transform~$\widehat\chi_B$.
For this, one needs to approximate $\chi_B$ by a function in
a class where the Poisson formula applies, say either the Schwartz
class, at least a continuous integrable function whose Fourier transform
is integrable.
The value at~$0$ of the Fourier transform~$\widehat\chi_B$
is the volume~$V(B)$ of~$B$. If one can bound the remaining terms,
one can hope to give an upper bound for the difference $N(B)-V(B)$.
This is at least a way to attack the circle problem,
which leads to non-trivial upper bounds.

In the next Section,
we shall give more details on this spectral theoretic
approach to Manin's conjecture  in the case of compactifications 
of algebraic groups.
As we have indicated, it relies on the  generalization 
of the theory of Fourier series and Fourier integrals 
and of the Poisson formula, which concern the groups~$\R$ and~$\R/\Z$,
to general locally compact topological groups.
We therefore call it the Fourier--Poisson method.

\subsection{Some details on the Fourier--Poisson method}\label{sec.fourier}

\subsubsection{Algebraic groups and their compactifications}
Let $X$ be an \emph{equivariant compactification} of an algebraic group~$G$
over a number field~$F$.
This means that $X$ is a projective variety
containing $G$ as a dense open subset and that the left and right actions of~$G$
on itself extend to  actions on~$X$.

Then, $G(F)$ is a discrete subgroup of the adelic group~$G(\AD_F)$.
Moreover, for line bundles $L$ on~$X$, the height $H_L$ can often 
be extended to a function on $G(\AD_F)$: there are functions $H_{L,v}$
such that, for $\mathbf g\in G(F)\subset G(\AD_F)$,
\[ H_L(\mathbf g)=\prod_{v\in \Val(F)} H_{L,v}(g_v). \]
If we use the right hand side 
of this formula to extend~$H_L$ to a function on~$G(\AD_F)$,
we obtain that the height zeta function~$Z_G(L;s)$
is an average over the discrete group~$G(F)$
of a function $H_L^{-s}$ on~$G(\AD_F)$ and one can use
\emph{harmonic analysis} and the spectral decomposition
to understand~$Z_G(L;s)$.

Important cases of groups 
in which mathematicians have been able to use harmonic analysis
to prove the conjectures of \textsc{Batyrev--Manin--Peyre} are the following.
\begin{itemize}
\item \emph{algebraic tori}: then $X$ is called a \emph{toric variety},
see \cite{batyrev-t96};
\item \emph{vector groups}, \cite{chambert-loir-t2002};
\item the Heisenberg group of upper triangular matrices,
\cite{shalika-t2004}; 
\item reductive groups, embedded in the \emph{wonderful compactification}
defined by \textsc{De Concini--Procesi},
see \cite{shalika-tb-t2007}.
\end{itemize}

\subsubsection{The Poisson formula}
Commutative groups are of course simpler to study, because
irreducible representations are one-dimensional (\emph{i.e.},
are \emph{characters}). Let us fix some notation.
Let $\mu$ be a Haar measure on the locally compact group~$G(\AD_F)$, 
given as product of suitable Haar measures~$\mu_v$ on $G(F_v)$.
Let $\widehat{G_v}=\Hom(G(F_v),\mathbf U)$ be the group of characters of~$G(F_v)$
and 
$\widehat{G(\AD_F)}=\Hom(G(\AD_F),\mathbf U)$ be the group
of  characters of~$G(\AD_F)$, where $\mathbf U$ is the group of complex
numbers of modulus~$1$.
By restriction, any character~$\chi\in\widehat{G(\AD_F)}$ induces 
a character $\chi_v$ of $G(F_v)$, for each place~$v$.
Let $G(F)^\perp$ be the orthogonal of~$G(F)$,
that is the group of characters of~$G(\AD_F)$ which are trivial on~$G(F)$.
This  is a locally compact group and it carries a dual measure  $\mu^\perp$.

Any integrable function $f\in L^1(G(\AD_F))$ has a \emph{Fourier transform},
which is the function on $\widehat{G(\AD_F)}$ defined as follows:
for $\chi\in  \widehat{G(\AD_F)}$, 
\[ \mathscr F(f;\chi)=\int_{G(\AD_F)} f(\mathbf g)\chi(\mathbf g)\,\mathrm d\mu(\mathbf g)  .\]
If $f$ is a simple function, \emph{i.e.}, of the form $f(\mathbf g)=\prod f_v(g_v)$,
then $\mathscr F(f;\chi)$ is a product 
$\prod\mathscr F_v(f_v;\chi_v)$ of local Fourier transforms
defined in an analogous manner, namely
\[ \mathscr F_v(f_v;\chi_v)= \int_{G(F_v)} f_v(g)\chi_v(g)\,\mathrm d\mu_v(g).\]

The generalization of \emph{Poisson formula} to this context
states that  if $f$ is an integrable function on~$G(\AD_F)$
such that moreover $\mathscr F(f;\cdot) $ is integrable
on~$G(F)^\perp$, then
\[ \sum_{g\in G(F)} f(g) = \int_{G(F)^\perp} \mathscr F(f;\chi) \,\mathrm d\mu^\perp(\chi).\]

This gives a formula for the height zeta function
\[ Z_G(L;s) = \int_{G(F)^\perp} \mathscr F(H_L^{-s};\chi)\,\mathrm d\mu^\perp(\chi), \]
which in many cases proved to be a key tool towards establishing 
the desired analytic properties
of the height zeta function.

\subsubsection{The Fourier transform at the trivial character}
Let us consider the Fourier transform of $f=H_L^{-s}$  at
the trivial character, in other words, the function
\[ \Phi(s)= \mathscr F(H_L^{-s},\mathbf 1)= \int_{G(\AD_F)} H_L(g)^{-s}\,\mathrm d\mu(g). \]
The decomposition of the function~$H_L$ as a product of function
on the spaces~$G(F_v)$ implies a similar decomposition of~$\Phi$
as a the product $\Phi(s)=\prod_v \Phi_v(s)$, where
\[ \Phi_v(s) = \mathscr F_v(H_L^{-s},\mathbf 1)  
= \int_{G(F_v)} H_L(g)^{-s}\,d\mu_v(g). \]

In the above mentioned cases, these functions have been often computed
explicitly, resorting to properties of algebraic groups.
However, this can also be done in a very geometric manner thanks
to the \emph{fundamental observation}: 
the integral $\Phi_v(s)$ on $G(F_v)$ can be
viewed as a geometric analogue on $X(F_v)$
of \textsc{Igusa}'s local zeta functions.
Assume for simplicity that $L=\omega_X^{-1}$. 
In Section~\ref{sec.igusa} we will explain how to prove
the following facts (see \S\,\ref{subsubsec.conclusion}):
\begin{itemize}
\item $\Phi_v(s)$ converges for $\Re(s)> 0$;
\item $\Phi_v(s)$ has a meromorphic continuation to~$\C$;
\item for almost all places~$v$, $\Phi_v(s)$ can be explicitly
computed (\textsc{Denef}'s formula);
\item the product $\prod \Phi_v(s)$ converges for $\Re(s)>1$;
\item this product has a meromorphic continuation which is governed by
some Artin L-function.
\end{itemize}
We shall in fact explain that these facts apply in much more
generality than that of algebraic groups.

\subsubsection{Conclusion of the proof}
To prove \textsc{Manin}'s conjecture for an equivariant compactification
of an algebraic group using the Fourier method,
a lot of work still remains to be done.
One first has to establish a similar meromorphic continuation 
for other characters
(other representations in the non-commutative cases) and, then,
to integrate all these contributions over $G(F)^\perp$.
Of course, the last step requires that one obtains good upper bounds
in the first step, so that integration is at all possible.
Details depend however very much on the specific groups involved
and cannot be described here.

\section{Heights and Igusa zeta functions}\label{sec.igusa}

This final section aims at explaining the \emph{geometric computation}
of the Fourier transform of the height function via the
theory of \textsc{Igusa} zeta functions.
It essentially borrows from the introduction of our 
recent article~\cite{chambert-loir-tschinkel2008}.

\subsection{Heights and measures on adelic spaces}
\subsubsection{Heights}
Let $X$ be a projective variety over a number field~$F$.
Let $L$ be an effective divisor on~$X$, and let us endow
the associated line bundle~$\mathscr O_X(L)$ with an adelic metric.
This line bundle possesses a canonical global section~$\mathsec f_L$ 
whose divisor is~$L$. For any place~$v\in \Val(F)$,
the function $x\mapsto \norm{\mathsec f_L(x)}_v$
is continuous on~$X(F_v)$ and vanishes precisely on~$L$.
Moreover, its sup-norm is equal to~$1$ for almost all places~$v$.

Consequently, the infinite product 
$\prod_{v\in\Val(F)} \norm{\mathsec f_L(x_v)}_v$
converges to an element of~$\R_+$ for any adelic point $x=(x_v)\in X(\AD_F)$.
Let $U=X\setminus L$. If $x$ is an adelic point of~$U$,
then $\norm{\mathsec f_L(x_v)}_v\neq 0$ for all~$v$,
and is in fact equal to~$1$ for almost all~$v$.
We thus can define the \emph{height} of a point~$x\in U(\AD_F)$ to be equal to
\[ H_L((x_v)) = \prod_{v\in\Val(F)} \norm{\mathsec f_L(x_v)}_v^{-1}. \]
The resulting function on~$U(\AD_F)$ is continuous and admits
a positive lower bound.
Moreover, for any~$B>0$, the set of points~$x\in U(\AD_F)$
such that $H_L(x)\leq B$ is compact.

\subsubsection{Local measures (II)}
The measure~$\tau_{X,v}$ on~$X(F_v)$ 
defined in~\S\,\ref{subsubsec.local-measures.I}
gives finite volume to~$X(F_v)$ and the open set $U(F_v)$ has full measure.
We modify the construction as follows so as to define a Radon measure
on~$U(F_v)$ whose total mass 
(unless $L(F_v)=\emptyset$)
is now infinite:
\[ \mathrm d\tau_{(X,L),v}(x) = \frac{1}{\norm{\mathsec f_L(x)}_v} \mathrm d\tau_{X,v}(x). \]
For example, in the important case where $X$ 
is an equivariant compactification
of an algebraic group~$G$, $U=G$ and $-L$ is a canonical  divisor,
this measure is a Haar measure on the locally compact group~$G(F_v)$.

\subsubsection{Definition of a global Tamagawa measure}
The product of these local measures does not converge in general,
and to use them to construct a measure on the adelic space $U(\AD_F)$,
it is necessary to introduce convergence factors, \emph{i.e.},
a family $(\lambda_v)$ of positive real numbers such that
the product
\[ \prod_{\text{$v$ finite}} \lambda_v\tau_{(X,L),v}(U(\mathfrak o_v)) \]
converges absolutely. The limit will thus be a positive real number
(unless some factor $U(\mathfrak o_v)$ is empty,
which does not happen if we remove an adequate finite set of places
in the product). The existence of such factors is a mere
triviality, since one could take for $\lambda_v$
the inverse of $\tau_{(X,L),v}(U(\mathfrak o_v))$. 
Of course, we claim for a meaningful definition of a family $(\lambda_v)$,
based on geometric or arithmetic invariants of~$U$.

The condition under which we can produce such factors
is the following: \emph{$X$ is proper, smooth, geometrically connected, and 
the two cohomology groups $\mathrm H^1(X,\mathscr O_X)$
and $\mathrm H^2(X,\mathscr O_X)$ vanish.}
Let us assume that this holds and let us define 
\[ M^0= \mathrm H^0(U_{\bar F},\mathscr O^*)/\bar F^*
\quad\text{and}\quad M^1= \mathrm H^1(U_{\bar F},\mathscr O^*)/\text{torsion}; \]
in words, $M^0$ is the Abelian group of invertible functions on~$U_{\bar F}$
modulo constants, and $M^1$ is the Picard group of~$U_{\bar F}$
modulo torsion. By the very definition as a cohomology group
of~$U_{\bar F}$, they possess a canonical action of
the Galois group~$\Gamma_{F}$ of~$\bar F/F$.
Moreover, $M^0$ is a free $\Z$-module of finite rank. Indeed, 
to an invertible function on~$U_{\bar F}$, one may attach its divisor,
which is supported by the complementary subset, that is~$L_{\bar F}$.
This gives a morphism from $M^0$ 
to the free Abelian group generated by the irreducible
components of~$L_{\bar F}$ and this map is injective because $X$ is normal.
It follows  that $M^0$ is a free $\Z$-module of finite rank.
Moreover, under the vanishing assertion above, then $M^1$
is a free $\Z$-module of finite rank too.
We thus can consider the Artin L-functions of these two $\Gamma_{F}$-modules,
$\mathrm L(s,M^0)$ and $\mathrm L(s,M^1)$.

Using  \textsc{Weil}'s conjectures, proved by \cite{deligne74}, 
to estimate the volumes $U(\mathfrak o_v)$,
in a similar manner to what \textsc{Peyre} had  done 
to  define the global measure~$\tau_X$ on $X(\AD_F)$,
we prove that  the family $(\lambda_v)$ given by
\[  \lambda_v = \frac{\mathrm L_v(1,M^0)}{\mathrm L_v(1,M^1)} \]
if $v$ is a finite place of~$F$ and $\lambda_v=1$ if $v$ is archimedean
is a family of convergence factors. In other words, 
the infinite product of measures
\[\frac{\mathrm L^*(1,M^1)}{\mathrm L^*(1,M^0)} 
\prod_{v\in\Val(F)}  \frac{\mathrm L_v(1,M^0)}{\mathrm L_v(1,M^1)}\,\mathrm d\tau_{(X,L),v}(x_v) \]
defines a Radon measure on the locally compact space~$U(\AD_F)$.
We denote this measure by $\tau_{(X,L)}$ and call it the \emph{Tamagawa measure}
on~$U(\AD_F)$.

\subsubsection{Examples for the Tamagawa measure}
This construction coincides with \textsc{Peyre}'s construction if $L=\emptyset$,
since then $M^0=0$ (there are no non-constant invertible functions
on~$X$) and $M^1=\Pic(X_{\bar F})/\text{torsion}$.

It also coincides with the classical constructions for algebraic groups.
Let us indeed assume that $X$ is an equivariant compactification of
an algebraic group~$G$ and that $-L$ is a canonical divisor,
with $G=U=X\setminus L$. We have mentionned that $\tau_{(X,L),v}$
is a Haar measure on~$G(F_v)$. 
If $G$ is semi-simple, or if $G$ is nilpotent,
then $M^0=M^1=0$; in that case $\lambda_v=1$,
according to the fact that the global measure does not need any convergence
factor. 
However, if $G$ is a torus, then $M^0$ is the group of characters of~$G_{\bar F}$
and $M^1=0$; our definition coincides with \textsc{Ono}'s.
In all of these cases, the global measure we construct
is of course a Haar measure on the locally compact group~$G(\AD_F)$.

It also recovers some cases of homogeneous spaces $G/H$,
where $G$ and~$H$ are semisimple groups over~$\Q$,
under the assumption that $H$ has no nontrivial characters,
like those studied by~\cite{borovoi-r1995}.
In fact, the definition of Tamagawa measures on
such homogeneous spaces does not need any convergence factors,
as can be seen by the fact that $M^0=0$ (since invertible functions
on~$G$ are already constants) and $M^1=0$ (because the Picard
group of $G/H$ is given by characters of~$H$).

As a last example, let us recall that \cite{salberger98} had shown that 
the universal torsors over toric varieties do not need any convergence
factors, because the convergence factor of the fibers (a torus)
compensates the one of the base (a projective toric variety).
In our approach, this is accounted by the fact,
discovered by~\cite{colliot-thelene-sansuc1987}, that such torsors
have neither non-constant invertible functions, nor non-trivial line bundles,
hence $M^0=M^1=0$ and $\lambda_v=1$.

\subsection{An adelic-geometric analogue of Igusa's local zeta functions}
As we explained above, when using the Fourier method to elucidate
the number of points of bounded height on some algebraic varieties,
we are lead to establish the meromorphic continuation of the
function of a complex variable:
\[ Z\colon s\mapsto \int_{U(\A_F)} H_L(x)^{-s}\,\mathrm d\tau_{(X,L)}(x). \]
(This function was denoted~$\Phi$ in~\S\ref{sec.fourier}.)
By definition,
the measure~$\tau_{(X,L)}$ is (up to the convergence factors
and the global factor coming from Artin L-functions) the product of 
the local measures $\tau_{(X,L),v}$ on the analytic manifolds $U(F_v)$.
Similarly, the integrated function
$x\mapsto H_L(x)^{-s}$ is the product of the local functions
$x\mapsto \norm{\mathsec f_L(x)}_v^s$.
Consequently, provided it converges absolutely, this adelic integral
decomposes as a product of integrals over local fields, 
\[ Z(s)= \frac{\mathrm L^*(1,M^1)}{\mathrm L^*(1,M^0)}
   \prod_{v\in\Val(F)} \frac{\mathrm L_v(1,M^0)}{\mathrm L_v(1,M^1)} Z_v(s)  ,
\qquad Z_v(s) = \int_{U(F_v)} \norm{\mathsec f_L(x)}_v^s\,\mathrm d\tau_{(X,L),v}(x). \]

\subsubsection{Geometric Igusa zeta functions}
\textsc{Igusa}'s theory of local zeta functions studies the analytic properties
of integrals of the form 
\[ I_\Phi(s) = \int_{F_v^n} \Phi(x) \abs{f(x)}_v^{s-1}\, \mathrm dx, \]
where $f$ is a polynomial and $\Phi$ is a smooth function
with compact support on the affine $n$-space over~$F_v$.
(See~\cite{igusa2000} for a very good survey on this theory.)

Our integrals $Z_v(s)$ are straightforward generalizations of 
\textsc{Igusa}'s local zeta functions.
Indeed, since $\tau_{(X,L),v}=\norm{\mathsec f_L}_v^{-1}\tau_{X,v}$,
we have
\[ Z_v (s) =   \int_{U(F_v)} \norm{\mathsec f_L(x)}_v^{s-1}\,\mathrm d\tau_{X,v}(x)
=   \int_{X(F_v)} \norm{\mathsec f_L(x)}_v^{s-1}\,\mathrm d\tau_{X,v}(x)
 \]
because $L(F_v)$ has measure~$0$ in~$X(F_v)$.
In our case, we integrate on a compact analytic manifold, rather than a
test function on the affine space; as always in Differential
geometry, partitions of unity and local charts convert a global
integral into a finite sum of local ones.
The absolute value of a polynomial has been replaced by
the norm of a global section of a line bundle; again, in local
coordinates, the function~$\norm{\mathsec f_L}_v$ 
is of the form $\abs f_v \phi$, where $f$ is a local equation of~$L$
and $\phi$ a continuous non-vanishing function.
However, because of these slight differences, we call such integrals
\emph{geometric Igusa zeta functions}.

From now on, 
we will assume in this text that over the algebraic closure~$\bar F$,
$L$ is a divisor with simple normal crossings. This means that
the irreducible components of the divisor~$L_{\bar F}$ are smooth
and meet transversally; in particular, any intersection
of part of these irreducible components is either empty,
or a smooth subvariety of the expected codimension.
This assumption is only here for the convenience  of the computation. 
Indeed, by the theorem on resolution of singularities of~\cite{hironaka1965},
there exists a proper birational morphism $\pi\colon Y\ra X$,
with $Y$ smooth, which is an isomorphism on~$U$
and such that the inverse image of~$L$ (as a Cartier divisor)
satisfies this assumption.
We may replace~$(X,L)$ by $(Y,\pi^*L)$ without altering
the definitions of~$Z$, $Z_v$, etc..
Let $\mathscr A$ be the set of irreducible components of~$L$;
for $\alpha\in\mathscr A$, let $L_\alpha$ be the corresponding
component, and $d_\alpha$ its multiplicity in~$L$.
We thus have $L=\sum_{\alpha\in\mathscr A} d_\alpha L_\alpha$.

This induces a stratification of~$X$ 
indexed by subsets~$A\subset\mathscr A$, the stratum $X_A$ being given by
\[ X_A= \{x\in X\,;\, x\in L_\alpha \Leftrightarrow \alpha\in A\}. \]

In the sequel, we shall moreover assume that for any $\alpha\in \mathscr
A$, $L_\alpha$ is geometrically irreducible.
This does not cover all cases, but the general case can be
treated using a similar analysis, the Dedekind zeta function of~$F$ being
replaced by Dedekind zeta functions of finite extensions~$F_\alpha$
defined by the components~$L_\alpha$; 
I refer the interested reader to~\cite{chambert-loir-tschinkel2008} for details.

\subsubsection{Local analysis}
Let $x$ be a point in $X_A(F_v)$.
The ``normal crossings'' assumption on~$L$ implies that
on a small enough analytic neighbourhood~$\Omega_x$ of~$x$,
one may find local equations $x_\alpha$ of $L_\alpha$, for $\alpha\in A$,
and complete these local equations into a system of local coordinates
$((x_\alpha)_{\alpha\in A}; y_1,\dots,y_s)$, so that $s+\Card(A)=\dim X$.
In this neighbourhood of~$x$, the measure~$\tau$ can be written
$\omega(x;y)\,\mathrm dx\,\mathrm dy$, with $\omega$ a positive 
continuous function, and there is a continous and positive function
$\phi$ such that 
\[ \norm{\mathsec f_L}_v = \phi \cdot \prod_{\alpha\in A} \abs{x_\alpha}_v^{d_\alpha}
  \]
on~$\Omega_x$.
If $\Phi_x$ is a continuous function with compact support on~$\Omega_x$,
we thus have
\[ \int_{X(F_v)} \Phi_x \norm{\mathsec f_L}^{s-1}\,\tau_{(X,L),v}
= \int_{\Omega_x} (\phi^{s-1} \omega \Phi_x) \prod_{\alpha\in A}\abs{x_\alpha}_v^{(s-1)d_\alpha}\, \prod_{\alpha\in A}\mathrm dx_\alpha\, \mathrm dy_1\dots\mathrm dy_s. \]
By comparison with the integral
\[ \int_{F_v} \Phi(u) \abs{u}_v^s \, \mathrm du \]
which converges for $\Re(s)>-1$ if $\Phi$ is compactly supported
on~$F_v$, we conclude that 
\[ \int_{X(F_v)} \Phi_x \norm{\mathsec f_L}^{s-1}\,\tau_{(X,L),v}\]
converges absolutely for $\Re(s)>0$.

Introducing a partition of unity, we conclude that 
for any place $v\in\Val(F)$,
the integral $Z_v(s)$ converges for $\Re(s)>0$ and defines
a holomorphic function on this half-plane.

\subsubsection{Denef's formula}
Our understanding of the infinite product of the functions $Z_v(s)$
relies in the explicit computation of almost all of them.
The formula below is a straightforward generalization of a formula 
which \cite{denef87} used to prove that, when 
$f$ is a polynomial with coefficients in the number field~$F$
and $\Phi$ is the characteristic function of the unit polydisk in~$F_v^n$,
then \textsc{Igusa}'s local zeta function  $I_\Phi(s)$ is a rational function
of $q_v^{-s}$ whose degree is bounded when $v$ varies within
the set of finite places of~$F$.

Assuming that the whole situation has ``good reduction'' at 
a place~$v$, we obtain
the following formula,
where $k_v$ is the residue field of~$F$ at~$v$ and $q_v$
is its cardinality:
\[ Z_v(s) = \sum_{A\subset\mathscr A} 
         q_v^{-\dim X} 
        \Card(X_A(k_v)). 
           \prod_{\alpha\in  A} 
           \frac{q_v-1}{q_v^{1+d_\alpha(s-1)}-1}.
\]
As in~\cite{denef87}, this is a consequence of the fact
that for each point~$\tilde x\in X(k_v)$,
in the local analytic description of the previous
paragraph, we may find local coordinates~$x_\alpha$ and~$y$
such that $\phi\equiv \omega=1$ which parametrize
the open unit polydisk with ``center~$\tilde x$.''
That this is at all possible is more or less what is meant
by ``having good reduction''.

\subsubsection{Meromorphic continuation}
Let us analyse the various terms of this formula.
The stratum $X_A$ has codimension~$\Card(A)$ in~$X$ (or is empty);
this implies easily that
\[ q_v^{-\dim X} \Card(X_A(k_v)) = \mathrm O(q_v^{-\Card(A)}). \]
More precisely, since $L_\alpha$ is geometrically irreducible
and $X_\alpha$ is an open subset of it, it follows
from \cite{deligne74} that
\[ q_v^{-\dim X} \Card(X_{\alpha}(k_v)) = q_v^{-1} + \mathrm O(q_v^{-3/2}). \]
(In fact, the estimates of~\cite{lang-w54} suffice for that.)
Similarly, $X_\emptyset=X\setminus L=U$, and
\[ q_v^{-\dim X}\Card(X_\emptyset(k_v)) = 1+ \mathrm O(q_v^{-1/2}). \]
In fact, the vanishing assumption on the cohomology groups
$\mathrm H^1(X,\mathscr O_X)$ and $\mathrm H^2(X,\mathscr O_X)$
implies that the remainder is even $\mathrm O(q_v^{-1})$.
By an integration formula of~\cite{weil82},
this expression is precisely equal to
$\tau_{(X,L),v}(U(\mathfrak o_v))$.

We therefore obtain that for any $\eps>0$, there exists $\delta>0$
such that if $\Re(s-1)>-\delta$, then
\[ Z_v(s) = \tau_{(X,L),v}(U(\mathfrak o_v))
 \prod_{\alpha\in\mathscr A} (1-q_v^{-1-d_\alpha(s-1)})^{-1}
 \left(1+\mathrm O(q_v^{-1-\eps})\right) . \]
If we multiply this estimate by the convergence factors $\lambda_v$ that
we have introduced, and which satisfy $\lambda_v=1+\mathrm O(q_v^{-1})$,
we obtain that the infinite product
\[ \prod_{\text{$v$ finite}} \lambda_v Z_v(s) \prod_{\alpha\in\mathscr A}
  \zeta_{F,v}(1+d_\alpha(s-1))^{-1} \]
converges absolutely and uniformly for $\Re(s-1)>-\delta/2$; it defines
a holomorphic function $\Phi(s)$ on that half-plane.  We
have denoted by $\zeta_{F,v}$ the local factor at~$v$ of the Dedekind
zeta function of the number field~$F$.
Then, the equality
\begin{align*}
 Z(s)  & = \frac{\mathrm L^*(1,M^1)}{\mathrm L^*(1,M^0)}
     \prod_{v\in\Val(F)} \lambda_v Z_v(s)  \\
&  = \frac{\mathrm L^*(1,M^1)}{\mathrm L^*(1,M^0)} 
 \left( \prod_{\alpha\in\mathscr A} \zeta_F(1+d_\alpha(s-1))\right)
  \Phi(s) \prod_{\text{$v$ archimedean}} Z_v(s) \end{align*}
shows that $Z$ converges absolutely, and defines a holomorphic
function on $\{\Re(s)>1\}$.
Since the Dedekind zeta function~$\zeta_F$ has a pole of order~$1$
at~$s=1$, $Z(s)$ admits a meromorphic continuation
on the half-plane $\{\Re(s)>1-\delta/2\}$,
whose only pole is at $s=1$, with multiplicity $\Card(\mathscr A)$.

\subsubsection{The leading term}
Moreover,
\[ \lim_{s\ra 1}  (s-1)^{\Card(\mathscr A)}Z(s)
 = \frac{\mathrm L^*(1,M^1)}{\mathrm L^*(1,M^0)}
 \left( \prod_{\alpha\in\mathscr A} d_\alpha^{-1}\right)
  \zeta_F^*(1)^{\Card(\mathscr A)} 
  \Phi(1) \prod_{\text{$v$ archimedean}} Z_v(1) .\]
By definition of~$\Phi$, one has
\begin{align*}
\Phi(1) &= \prod_{\text{$v$ finite}} \lambda_v \prod_{\alpha\in\mathscr A}\zeta_{F,v}(1)^{-1} Z_v(1) \\
&=  \prod_{\text{$v$ finite}}
   \left( \lambda_v \prod_{\alpha\in\mathscr A}\zeta_{F,v}(1)^{-1} \right)
  \tau_{X,v}(X(F_v)).
\end{align*}
Comparing the Galois modules $M^0$ and~$M^1$ for $X$ and~$U$,
one can conclude that
\begin{equation}
\label{eq.Z1}
\lim_{s\ra 1}  (s-1)^{\Card(\mathscr A)}Z(s)
     = 
   \tau_X(X(\AD_F))
\prod_{\alpha\in\mathscr A} d_\alpha^{-1}.
\end{equation}
 
\subsubsection{Conclusion of the computation}\label{subsubsec.conclusion}
Finally, we have shown the following: The integral
\[ s\mapsto \int_{U(\AD_F)} H_L(x)^{-s}\,\mathrm d\tau_{(X,L)}(x) \]
converges for $\Re(s)>1$, defines a holomorphic function on this
half-plane. It has a meromorphic continuation
on some half-plane $\{\Re(s)>1-\delta\}$ (for some $\delta>0$)
whose only pole is at $s=1$, has order~$\Card(\mathscr A)$
and its asymptotic behaviour at $s=1$ is given by Equation~\eqref{eq.Z1}.

\subsection{Application to volume estimates}
Let us finally derive our volume estimates from
the analytic property of the function $Z(s)$.
We are interested in the volume function defined by
\[ V(B) = \tau_{(X,L)}\left(\{ x\in U(\AD_F)\,;\, H_L(x)\leq B\}\right) \]
for $B>0$, more precisely in its asymptotic behaviour when
$B\ra\infty$.

Its Mellin-Stieltjes transform is given by
\[  \int_0^\infty B^{-s} \mathrm dV(B) 
= \int_{U(\AD_F)}  H_L(x)^{-s}\,\mathrm d\tau_{(X,L)}(x) = Z(s). \]
By a slight generalization of \textsc{Ikehara}'s Tauberian theorem,
we conclude that $V(B)$ satisfies the following asymptotic expansion:
\[ V(B) \sim  \big((a-1)!\prod_{\alpha\in\mathscr A} d_\alpha \big)^{-1}
  \tau_X(X(\AD_F)) B (\log B)^{a-1}, \]
where $a=\Card(\mathscr A)$ is the number of irreducible
components of the divisor~$L$.

Our argument, applied to different metrizations, also implies
that when $B\ra\infty$, the probability measure
\[ \frac 1{V(B)} \mathbf 1_{H_L(x)\leq B} \mathrm d\tau_{(X,L)}(x) \]
on $U(\AD_F)$ (viewed as a subset of $X(\AD_F)$)
converges to the measure
\[ \frac1{\tau_{X}(X(\AD_F))} \,\mathrm d\tau_X(x) \]
on $X(\AD_F)$.

\bigskip
\subsubsection*{Acknowledgments}
I would like to thank Yuri \textsc{Tschinkel}
for his comments on a first version of this survey.
Our collaboration on the topics presented here
began more than ten years ago.  Let it be the occasion 
to thank him heartily 
for having shared his views and projects with me.

I also thank the referee, as well as R.~\textsc{de la Bretèche},
H.~\textsc{Oh} and W.~\textsc{Veys} for their suggestions.

\bibliographystyle{mynat}
\bibliography{aclab,acl,miura}

\end{document}